\documentclass[12pt]{article}

\usepackage{amsmath,amssymb,amsthm}
\usepackage{a4wide}
\usepackage[usenames,dvipsnames]{xcolor}
\usepackage{graphicx}
\usepackage{natbib}
\usepackage{mathrsfs}
\usepackage{bm}

\usepackage{mathrsfs}
\usepackage{verbatim}
\usepackage{slashbox}
\usepackage{slashbox}

\makeatletter
\newtheorem{lemma}{Lemma}[section]
\newtheorem{theorem}{Theorem}[section]
\newtheorem{remark}{Remark}
\newtheorem{example}{Example}[section]
\newtheorem{corollary}{Corollary}[section]

\numberwithin{equation}{section}
\allowdisplaybreaks       

\newcommand\mgai[1]{{#1}}

\newcommand\gai[1]{{#1}}  
\begin{document}

\def\di{D^{-1}}
\def\dij{D^{-1}_j}

\newcommand{\bqa}{\begin{eqnarray}}
  \newcommand{\eqa}{\end{eqnarray}}
\newcommand{\bqn}{\begin{equation}\begin{array}{lll}}
    \newcommand{\eqn}{\end{array}\end{equation}}
\newcommand{\um}{\underline{m}}
\newcommand{\be}{\begin{equation}}
  \newcommand{\ee}{\end{equation}}
\newcommand{\sjln}{{\sum\limits_{j=1}^{n}}}
\newcommand{\siln}{{\sum\limits_{i=1}^{n}}}
\newcommand{\bgma}{\bm{\gamma}}
\newcommand{\non}{\nonumber\\}
\newcommand{\bbA}{{\bf A}}
\newcommand{\bbC}{{\bf C}}
\newcommand{\bbD}{{\bf D}}
\newcommand{\rE}{{\rm E}}
\newcommand{\bbK}{{\bf K}}
\newcommand{\bbB}{{\bf B}}
\newcommand{\bbI}{{\bf I}}
\newcommand{\bbr}{{\bf r}}
\newcommand{\bbq}{{\bf q}}
\newcommand{\bbQ}{{\bf Q}}
\newcommand{\bbM}{{\bf M}}
\newcommand{\bbU}{{\bf U}}
\newcommand{\bbV}{{\bf V}}
\newcommand{\bgL}{\bm{\Lambda}}
\newcommand{\bbS}{{\bf S}}
\newcommand{\bbT}{{\bf T}}
\newcommand{\bbx}{{\bf x}}
\newcommand{\bbX}{{\bf X}}
\newcommand{\bby}{{\bf y}}
\newcommand{\rtr}{{\rm tr}}
\newcommand{\bA}{{\bf A}}
\newcommand{\bT}{{\bf T}}
\newcommand{\bbgS}{{\boldsymbol\Sigma}}
\newcommand\re[1]{{\color{red}#1}}
\newcommand{\cC}{{\mathcal{C}}}
\newcommand{\ep}{{\epsilon}}
\newcommand{\rCov}{{\rm Cov}}
\newcommand{\Cov}{{\rm Cov}}
\newcommand{\rVar}{{\rm Var}}
\newcommand{\darrow}{\downarrow}
\newcommand{\iparrow}{\stackrel{i.p.}{\rightarrow}}
\newcommand{\umn}{\underline{m}_n}
\newcommand{\diag}{{\rm diag}}
\newcommand{\rank}{{\rm rank}}
\newcommand{\rP}{{\rm P}}
\newcommand{\gs}{\sigma}
\newcommand{\hSigma}{\hat{\bm{\Sigma}}}
\newcommand{\bSigma}{\bm{\Sigma}}
\renewcommand\hat[1]{\widehat{#1}}
\renewcommand\tilde[1]{\widetilde{#1}}

\title{CLT for linear spectral statistics of large dimensional sample covariance matrices with dependent data}
\author{Shurong Zheng, Zhidong Bai, Jianfeng Yao and Hongtu Zhu}
\date{July 13, 2016}
\maketitle{}
\begin{footnotetext}
  {Shurong Zheng and Zhidong Bai are professors in School of Mathematics $\&$ Statistics and KLAS, Northeast Normal University, China;
    Jianfeng Yao is professor \mgai{with the}  Department of Statistics and Actuarial Science, The University of Hong Kong,  Hong Kong;
    Hongtu Zhu is professor \mgai{at the} University of North Caroline at Chapel Hill, USA.}
\end{footnotetext}

\begin{abstract}
  This paper investigates the central limit theorem for linear spectral statistics of high dimensional sample covariance matrices of the form $\bbB_n=n^{-1}\sum_{j=1}^{n}\bbQ\bbx_j\bbx_j^{*}\bbQ^{*}$
  where $\bbQ$ is a nonrandom matrix of dimension  $p\times k$, and
  $\{\bbx_j\}$ is
\mgai{a sequence of independent}  $k$-dimensional random vector
\mgai{with} independent entries, under the assumption that $p/n\to
y>0$. A key \mgai{novelty here}
is that the dimension \mgai{$k\ge p$}  can be arbitrary, \mgai{possibly}  infinity.
This new model of sample covariance matrices $\bbB_n$ covers most of
the known models as its special cases.
\mgai{For example,  standard sample covariance matrices are obtained
  with  $k=p$ and $\bbQ=\bbT_n^{1/2}$  for some positive definite
  Hermitian matrix $\bbT_n$.  Also with $k=\infty$ our model covers
  the case of  repeated linear  processes considered in recent
  high-dimensional time series literature.  The CLT}
found in this paper
substantially generalizes \mgai{the seminal CLT} in  Bai and
Silverstein (2004).
\mgai{Applications of this new CLT are proposed for
  testing the structure of a high-dimensional covariance matrix.
  The derived tests are then used to analyse a
  large fMRI data set regarding its temporary correlation structure.}
\end{abstract}


\section{Introduction}
\label{sec:intro}
\mgai{Sample covariance matrices are of central importance in multivariate
  and high-dimensional data analysis.
  They have prominent applications in various big-data fields
  such as wireless communication networks, social networks and signal
  processing,
  see \mgai{for instance} the recent survey papers \cite{CM2013},
  \cite{John07} and \cite{PaulAue14}.
  Many statistical tools in these area
  depend on the so-called {\em linear spectral statistics} (LSS)  of a sample
  covariance matrix, say $\bbB_n$ of size $p\times p$ with eigenvalues $\{\lambda_i\}_{i=1}^{p}$,
  which have the form $p^{-1}\sum_{i=1}^{p}f(\lambda_i)$
  where $f$ is a given function.}
Two main questions arise for such  LSS, namely
(i)  determining the point limit of a {\em limiting spectral distribution} (LSD), say $G$, such that LSS converges to $G(f)=\int f(x) dG(x)$
(in an appropriate sense and for a wide family of functions $f$);
and  (ii) characterizing the fluctuations $p^{-1}\sum_{i=1}^{p}f(\lambda_i)-G(f)$ in terms of an appropriate central limit theorem (CLT).
Both questions  have a long history and \mgai{continue to receive}
considerable attention in recent years.
As for LSDs, the question has been extensively studied in the literature starting
from the classical work of   Mar$\check{c}$enko and Pastur (1967), and
contitued in \cite{Silv95} and \cite{Wachter1978}, and
found many recent developments
\citep[among others.]{BZ2008,yao2012,BM2015}

As for CLTs for linear spectral statistics,
a CLT for $(\rtr(\bbB_n),$ $\cdots,$ $\rtr(\bbB_n^{\ell}))$ is established  in  \cite{Jons82}
for a sequence of Wishart matrices $\{\bbB_n\}$ with Gaussian variables, where $\ell$ is a fixed number, and the dimension $p$ of the matrices grows
proportionally to the sample size $n$. By employing a general method based on Stieltjes transform,
\cite{BS04} provides a CLT for LSS of large dimensional sample covariance matrices from a general population,
that is, not  necessarily Gaussian  along  with an arbitrary population covariance matrix.
A distinguished feature of this CLT is that its centering term and
limiting means and covariance functions are  all fully
characterized.
\mgai{This celebrated paper has been subsequently} improved
in several follow-up papers including
\cite{PanZhou08}, \cite{Pan12}, \cite{NJ2013}  and \cite{ZBY2015}.
\mgai{These CLTs have found successful applications
  in high-dimensional statistics
  by  solving notably important problems in parameter estimation or
  hypothesis testing
  with high-dimensional data. Examples include}
\cite{bjyz09},  \cite{WY13} and \cite{lxzt14} for hypothesis testing on
covariance matrix;
\cite{bjyz13} for testing on regression coefficients,
and \cite{jbz13} for testing indepedence between two large sets of variables, among  others.

The aim of this paper is to establish a new CLT
for a \mgai{much extended}  class  of sample covariance matrices. This
new development  represents a  novel  and major extension of the   CLT
of \cite{BS04} \mgai{and its recent follow-up version}.
Specifically,
in this paper, we consider samples   $\{\bby_j\}_{1\le  j\le n}$ \mgai{of
the form}  $\bby_j=\bbQ\bbx_j$ where
\begin{enumerate}
\item[(M1)] $\{\bbx_j=(x_{1j},\ldots,x_{kj})^T,~1\le j\le n\}$ is a
  sequence of \mgai{independent} $k$-dimensional random
  vectors with independent standardized components $x_{ij}$,
  i.e. $Ex_{ij}=0$ and $E|x_{ij}|^2=1$, and the dimension $k\ge
    p$, {\em possibly  $k=\infty$}
    (for example, each $\bbx_j$ is a white noise \mgai{made with} an
    infinite sequence of standardized \mgai{errors});
\item[(M2)]  $\bbQ$ is a $p\times k$ matrix  \emph{with arbitrary entries}.
\end{enumerate}
The sample covariance matrix in our setting is then given by
\begin{equation}
  \label{eq:Bn}
  \bbB_n=n^{-1}\sum_{j=1}^{n}\bby_j \bby_j^*=n^{-1}\bbQ \bbX_n \bbX_n^*\bbQ^*~,
\end{equation}
where $\bbX_n=(\bbx_1,\ldots,\bbx_n)$ is an $k\times n $ matrix (with
independent entries) and $*$ stands for the transpose and complex
conjugate of matrices or vectors. In contrast,
\mgai{classical sample covariance matrices as in}
\cite{BS04} consider samples of the form  $\bby_j={\bf C}_n^{1/2}\bbx_j$,  where
\begin{enumerate}
\item[(O1)] $\{\bbx_j=(x_{1j},\ldots,x_{pj})^T,~1\le j\le n\}$ is a sequence of iid $p$-dimensional random
  vectors with independent standardized components $x_{ij}$, i.e. $Ex_{ij}=0$ and $E|x_{ij}|^2=1$.
\item[(O2)] ${\bf C}_n^{1/2}$ is a $ p\times p$ positive definite Hermitian
  matrix.
\end{enumerate}
The sample covariance matrix of interest thus takes form
\begin{equation}
  \label{eq:Sn}
  \bbS_n=n^{-1}\sum_{i=1}^n\bby_j\bby_j^*=n^{-1}{\bf C}_n^{1/2}\bbX_n\bbX_n^*{\bf C}_n^{1/2},
\end{equation}
where $\bbX_n=(\bbx_1,\ldots,\bbx_n)$ is a $p\times n$ data matrix (with independent entries).

\gai{The main innovation} \mgai{in the model} \eqref{eq:Bn} is that the {\em mixing} matrix $\bbQ$ can have an infinite number of columns, and this will allow to
cover dependent samples. For example, the repeated linear process
$\{y_{ij}=\sum_{t=-\infty}^{\infty}b_tx_{i-t,j}, i=1,\ldots,p,
j=1,\ldots,n\}$ is a special case
of  (\ref{eq:Bn}) with
\[
\bby_j= \bbQ\bbx_j= \underbrace{\left(\begin{array}{cccccccc}
      \cdots&b_{1-p}&b_{2-p}&\ldots&b_0&b_1&b_2&\ldots\\
      \vdots&\vdots&\vdots&\vdots&\vdots&\vdots&\vdots\\
      \cdots&b_{-1}&b_0&\ldots&b_{p-2}&b_{p-1}&b_{p}&\ldots\\
      \cdots&b_0&b_1&\ldots&b_{p-1}&b_{p}&b_{p+1}&\ldots\\
    \end{array}\right)}_{\bbQ~\mbox{is}~p\times \infty~\mbox{dimensional}}
\underbrace{\left(\begin{array}{c}
      \vdots\\
      x_{p,j}\\
      \vdots\\
      x_{1,j}\\
      \vdots
    \end{array}\right)}_{\bbx_j~\mbox{is}~\infty-\mbox{dimensional}}.
\]
Such data structure often arises in panel surveys or longitudinal studies, in which  respondents or subjects are   
interviewed/observed during a certain time period (see \cite{WongMiller1990, WMS2001}).
Moreover, if indeed $k=p$, then (\ref{eq:Sn}) of \cite{BS04} is a special case of our novel extension. In this case,  although   the entries of $p\times p$ matrix $\bbQ$  are arbitrary and $\bbQ$ is not required to be Hermitian or positive definite,
our setup is in fact strictly equivalent to \cite{BS04}
due to  the  polar decomposition  $\bbQ=\bbU_n\bbT_n^{1/2}$, in which  $\bbU_n$ is
unitary and  $\bbT_n^{1/2}$  is  Hermitian and nonnegative definite.
Consequently, $\bbB_n=\bbU_n{\bbS}_n\bbU_n^*$
where $\bbS_n$ is equivalent to the sample covariance matrix
in \eqref{eq:Sn} has been considered in \cite{BS04}.
 Therefore, $\bbB_n$ and $\bbS_n$ have exactly the same spectrum of
eigenvalues and their LSS are identical. However,  when $k$ tends to
infinity with a higher order than $p$,
 the new  model in \eqref{eq:Bn} represents
a  non-trivial extension of the framework of \cite{BS04}.

We will derive  the LSD of the sample  covariance matrix $\bbB_n$
in \eqref{eq:Bn} and the corresponding CLT for its linear spectral
statistics, thus extending  both the results of \cite{BS04} and
\cite{Silv95} \mgai{into the new framework}.
As it will be seen below,
establishing  these extensions  requires
several novel  ideas and major techniques
\mgai{in order to tackle with a}
mixing matrix $\bbQ$ that can be infinite and with arbitrary  entries.

The rest of the paper is organized as follows.
In Section~\ref{sec:LSD}, we derive the LSD  of the sample covariance matrix
$\bbB_n$ under suitable conditions.
In Section~\ref{sec:CLT}, we derive the CLT for linear spectral statistics of
$\bbB_n$ under some reinforced conditions.
In Section~\ref{sec:app}, we present two applications   of these theoretical results.
\mgai{Section~\ref{sec:proofs} collects the proofs for the main theorems of
the paper.}
Finally, some  technical  lemmas and auxiliary results used in the
proofs
of   Section~\ref{sec:proofs}  are postponed to the
appendices.

\section{LSD of the sample covariance matrix $\bbB_n$}
\label{sec:LSD}

This section aims at the derivation of the
LSD  for the sample covariance matrix
$\bbB_n$ in \eqref{eq:Bn} from sample
$\bby_j$'s of dependent data.
Recall that the {\em empirical spectral distribution} (ESD)
of a $p\times p$ square  matrix $\bbA$ is the probability
measure $F^\bbA=p^{-1} \sum_{i=1}^p\delta_{\lambda_i}$,
where the $\lambda_i$'s  are eigenvalues of $\bbA$ and $\delta_a$ denotes the
Dirac mass at point $a$. For any probability measure $F$ \mgai{on the real line},
its Stieltjes transform is defined by
\[
m(z)=\int\frac{1}{t-z} dF(x), \quad  z\in
\mathbb{C}^+=\{z\in\mathbb{C}:~ \Im(z)>0\},
\]
where \mgai{$\mathbb{C}$  denotes the complex plane and $\Im$ the imaginary part}.

The assumptions needed for the derivation of the LSD of $\bbB_n$ are as follows.
\begin{description}
\item[Assumption (a)]
  Samples are $\{\bby_j=\bbQ\bbx_j, j=1,\ldots,n\}$,  where $\bbQ$ is $p\times k$, $\bbx_j$ is $k\times 1$, $\bbx_j=(x_{1j},\ldots,x_{kj})^T$,
  and  \emph{the dimension $k$ is arbitrary (possibly infinite)}. Moreover,
  $\{x_{ij}, i=1,\ldots,k, j=1,\ldots,n\}$ is a $k\times n$ array of
  independent random variables, not necessarily
  identically distributed, with common moments
  $$
  E x_{ij}=0, ~~ E|x_{ij}^2|=1, 
  $$
  and satisfying the following Lindeberg-type condition:
  for each $\eta>0$,
  $$
  \frac1{pn\eta^2}\sum_{i=1}^k\sum_{j=1}^n \|\bbq_{i}\|^2\rE|x_{ij}^2|I\Big(|x_{ij}|>\eta \sqrt{n}/\|\bbq_i\|\Big)\to0,
  $$
  where $\|\bbq_i\|$ is the Euclidean norm
  \mgai{($k$ is finite)  or the $\ell^2$ norm ($k=\infty$)} of the $i$-th column vector
  $\bbq_i$ of $\bbQ$.
\item[Assumption (b)] With probability one, the ESD  $H_n$ of the
  population covariance matrix $\bbT_n=\bbQ\bbQ^{*}$ converges weakly
  to a probability distribution $H$.
  \mgai{Also the sequence $(\bbT_n)_n$ is  bounded in spectral  norm.}
\item[Assumption (c)] Both $p$ and $n$ tend to infinity such  that $y_n=p/n\to y>0$ as $n\to\infty$.
\end{description}

The Lindeberg-type condition in Assumption {\bf (a)} is a classical
moment condition of second order. It is automatically satisfied if
the entries $\{x_{ij}\}$ are identically distributed. This  condition
is here used to tackle with possibly non identically distributed
entries and ensures a suitable truncation of the variables.
Assumption {\bf (b)} is also a standard condition on the convergence
of the population spectral distribution and
requires the weak convergence of the ESD of the population covairance matrix.
Assumption {\bf (c)} defines the asymptotic regime following the
seminal work of Mar\v{c}enko and Pastur (1967).

As the  first main result of the paper, we derive the LSD of $\bbB_n$.

\begin{theorem}\label{thm1}
  Under Assumptions \textbf{(a), (b)} and \textbf{(c)}, almost surely the  ESD
  $F^{\bbB_n}$ of $\bbB_n$ weakly converges  to a non-random
  LSD  $F^{y,H}$. Moreover, the LSD $F^{y,H}$ is determined through its
  Stieltjes transform $m(z)$,  which is the unique
  solution to the following Mar\v{c}enko-Pastur equation
  \begin{equation}
    \label{eq:MP}
    m(z)= \int\frac{1}{t[1-y-yzm(z)]-z}dH(t)~,
  \end{equation}
  on the set
  $\{ m(z) \in\mathbb{C}:  -(1 - y)/z + ym(z)\in\mathbb{C}^+\}$.
\end{theorem}

Theorem~\ref{thm1}  has several important implications. The LSD $F^{y,H}$ is exactly the generalized Mar\v{c}enko-Pastur
distribution with parameters $(y,H)$  as described in \cite{Silv95}. The key  novelty of Theorem~\ref{thm1} consists of an
extension of this LSD into a much more general setting on    matrix entries as defined  in Assumptions {\bf (a)-(b)}. Therefore,  Theorem~\ref{thm1}    includes many existing results related to LSDs of sample covariance matrices from dependent data as special cases, e.g. \cite{BZ2008,JWBNH2014,wjm09,wjm11,yao2012}.
To the best of our knowledge,  although  the  model considered in  \cite{BZ2008} is more general than that in  this paper, it relies on  a strong quadratic form condition,  which can be hardly  verified in many applications.

Further developments on this LSD rely on an equivalent representation of the Mar\v{c}enko-Pastur law.
Namely, define the {\em companion LSD} of $\bbB_n$ as
$$
\underline{F}^{y,H}=(1-y)\delta_0+yF^{y,H}.
$$
It is readily checked that
$ \underline{F}^{y,H}$ is the LSD of the {\em companion sample
  covariance matrix} $\underline{\bbB}_n=n^{-1}\bbX_n^{*}\bbQ^{*}\bbQ\bbX_n$ (which is $n\times n$), and its Stieltjes transform
$\underline{m}(z)$ satisfies \mgai{the so-called Silverstein equation}
\begin{equation}
  z=-\frac1{\um(z)} +y\int\frac{t}{1+t\um(z)}dH(t), \quad \underline{m}(z)=-\frac{1-y}{z}+ym(z).
  \label{eq1}
\end{equation}
The advantage of this  equation is that it indeed defines the {\em inverse function} of the Stieltjes transform $\um(z)$.
This equation is the key to numerical evaluation of the Stieltjes transform, or the underlying density function of the
LSD.  Moreover, many analytical properties on LSD can be inferred from this equation (see \cite{SC1995}).

We consider an example to elaborate more on Theorem~\ref{thm1}.

\begin{example}
  Assume that for each $1\le j\le n$,  $\{y_{ij}\}_i$ is an ARMA(1,1) process of the  form
  \[
  y_{ij}=\phi y_{i-1,j}+\theta\varepsilon_{i-1,j}+\varepsilon_{ij},
  \]
  where $\{\varepsilon_{ij}\}$ is an array of independently and identically distributed   random variables with mean zero and
  variance 1. {The condition $|\phi|\vee |\theta|<1$ is assumed for causality and invertibility of this ARMA(1, 1) process.}
  From (\ref{eq:MP}) or (\ref{eq1}), it follows that the Stieltjes transform $m(z)$ of the LSD $F^{y,H}$  satisfies the following equation
  $$
  z=-\frac{1}{m(z)}+\frac{\theta}{y\theta m(z)-\phi}-\frac{(\phi+\theta)(1+\phi\theta)}{(y\theta m(z)-\phi)^2}\frac{\epsilon(\alpha)}{\sqrt{\alpha^2-4}}
  $$
  with
  $$
  \alpha=\frac{ym(z)(1+\theta^2)+1+\phi^2}{y\theta m(z)-\phi}~~\mbox{and}~~ \epsilon(\alpha)={\rm sgn}(\Im(\alpha)).
  $$
  In the case of an AR(1) process ($\theta=0$), we have
  $$
  z=-\frac{1}{m(z)}+\frac{1}{\sqrt{[y m(z)+1+\phi^2]^2-4\phi^2}}.
  $$
  Similarly, in the case of a MA(1) processs ($\phi=0$),  we have
  $$
  z=-\frac{1}{m(z)}+\frac{1}{ym(z)}+\frac{1}{y^2m^2(z)}\frac{1}{\sqrt{\{[ym(z)]^{-1}+1+\theta^2\}^2-4\theta^2}}
  $$
  (see \cite{wjm09}, \cite{wjm11} and \cite{yao2012}).
\end{example}

\begin{remark}
  In general,  we cannot find a closed-form formula for the Stieltjes transform $m(z)$ (or $\um(z)$). One has to rely on numerical evaluation
  using the Silverstein equation \eqref{eq1}. Notice that the values of $m(z)$ at points $z=x+{\bf i}\varepsilon$ with small positive
  $\varepsilon$ provides an approximation  to the value of the density of the LSD at $x$ by the inversion formula.
  For details on these numerical algorithms, we refer to
  \cite{Dobrian15}. By (\ref{eq1}) and  \cite{SC1995},  we obtain $\underline{m}(z)$ and   the limiting density $f^{y,H}$ of $F^{y,H}$, which  satisfies
  $f^{y,H}(x)=(y\pi)^{-1}\lim\limits_{z\rightarrow x+0{\bf i}}\Im(\underline{m}(z))$.
\end{remark}

\begin{remark}\label{rem001}
  If ${\bf T}_n=\bbQ\bbQ^{*}$ is the identity matrix $\bbI_p$, then
  \mgai{$H=\delta_1$ and we have by \eqref{eq1},}
  $$
  z=-\frac1{\um(z)}+\frac{y}{1+\um(z)}.
  $$
  That is,
  $$
  \um(z)=\frac{-(z+1-y)+\sqrt{(z-1-y)^2-4y}}{2z}
  $$
  and the LSD $F^{y,H}$ of $\bbB_n$ is the Mar$\check{c}$enko-Pastur law with the density function
  \begin{eqnarray}
    f^{y,H}=(2\pi yx)^{-1}\sqrt{[(1+\sqrt{y})^2-x][x-(1-\sqrt{y})^2]}
  \end{eqnarray}
  and has a point mass $1-1/y$ at the origin if $y>1$.
\end{remark}

\section{CLT for linear spectral statistics of the  sample covariance matrix $\bbB_n$}
\label{sec:CLT}
In this section,
we establish the corresponding CLT for LSS of
\mgai{the sample covariance matrix $\bbB_n$ in \eqref{eq:Bn} under approriate conditions.}
This CLT constitutes  the second main contribution of this paper.
\mgai{We first introduce the assumptions needed for this CLT}.
\begin{description}
\item[Assumption (d)]
  The variables $\{x_{ij}, i=1,\ldots,k, j=1,\ldots,n\}$ are independent, with common moments
  $$
  E x_{ij}=0, ~~ E|x_{ij}^2|=1, ~~
  \beta_x=E|x_{ij}^4|-|Ex_{ij}^2|^2-2, ~~\mbox{and}~~\alpha_x=|Ex_{ij}^2|^2,
  $$
  and satisfying the following Lindeberg-type condition: for each $\eta>0$
  {\be
    \frac1{pn\eta^6}\sum_{i=1}^k\sum_{j=1}^n\|\bbq_{i}\|^2\rE|x_{ij}^4|I\Big(|x_{ij}|>\eta\sqrt{n/\|\bbq_i\|}\Big)\to0.
    \ee}

\item[Assumption (e)]
  Either $\alpha_x=0$,  or  the mixing matrix
  $\bbQ$ is real (with arbitrary $\alpha_x$).
\item[Assumption (f)]
  {Either  $\beta_x=0$, or the mixing matrix $\bbQ$ is such that diagonal $\bbQ^*\bbQ$ (with arbitrary $\beta_x$).}
\end{description}

\begin{remark}
  Assumption {\bf (d)}  is  stronger than  Assumption {\bf (a)} in
  Section~\ref{sec:LSD} since Assumption {\bf (d)}  assumes  the
  Lindeberg-type condition on the
  fourth order moments (instead of the previous moments of second
  order). Assumption {\bf (e)} means that if variables $x_{ij}$s' are
  complex-valued,
  $E(x_{ij}^2)=0$ should hold which is the Gaussian-like second moment
  in \cite{BS04}. However,  if the other real $\bbQ$ condition is
  imposed, then the Gaussian-like second moment can be
  removed ($\alpha_x\ne 0$).
  \mgai{If $\alpha_x\ne 0$ and
    $\bbQ$ is not real, then the CLT for LSS of $\bbB_n$ may not
    hold}.
  Such  counterexample for the case of $k=p$ can be found
  in \cite{ZBY-cltF}. As for Assumption {\bf (f)}, if the population
  is Gaussian, then
  $\beta_x=0$  is the Gaussian-like fourth moment. If we impose the condition that the matrix $\bbQ^*\bbQ$ is real and
  diagonal, then the Gaussian-like fourth moment condition can be
  removed.
  \mgai{Again if $\beta_x\ne 0$ and $\bbQ^*\bbQ$ is not diagonal,
    the CLT may not hold.
  }
  Such  counterexamples with  $k=p$ can be found in \cite{ZBY-cltF}.
\end{remark}

\begin{theorem} \label{thm2}
  Under Assumptions ({\bf b})-({\bf f}), let  $f_1,\ldots,f_L$ be $L$ functions analytic on a complex domain containing
  \be
  [I_{(0<y<1)}(1-\sqrt{y})^2\liminf_n\lambda_{\min}^{\bbT_n},~~(1+\sqrt{y})^2\limsup_n\lambda_{\max}^{\bbT_n}]\label{tag1.4} \ee
  with $\bbT_n=\bbQ\bbQ^{*}$,  and $\lambda_{\min}^{\bbT_n}$ and
  $\lambda_{\max}^{\bbT_n}$ \mgai{denoting its smallest and the
    largest  eigenvalue, respectively}.
  Consider
  the random vector $\left(X_p(f_1),\ldots,X_p(f_L)\right)$ with
  \begin{equation}\label{Yp}
    X_p(f_\ell)=\sum_{j=1}^p f_\ell(\lambda_j) - pF^{y_{n},H_n}(f_\ell),\quad \ell=1,\ldots, L,
  \end{equation}
  in which $F^{y_{n},H_n}(f_\ell)=\int f_{\ell}(x)f^{y_n,H_n}(x)dx$ and $\{\lambda_j\}_{j=1}^{p}$ are the sample eigenvalues of $\bbB_n$.
  Then, the random vector $(X_p(f_1),\ldots,X_p(f_L))$ converges to a $L$-dimensional Gaussian random vector $(X_{f_1},\ldots,X_{f_L})$ with mean function
  \begin{eqnarray}
    \rE X_{f_\ell}&=&-\frac{1}{2\pi {\bf i}}\oint\limits_{\cal C}
    f_\ell(z)\frac{\alpha_xy\int\frac{\underline{m}^3(z)t^2}{(1+t\underline{m}(z))^{3}}dH(t)}
    {\left(1-y\int\frac{\underline{m}^2(z)t^2}{(1+t\underline{m}(z))^{2}}dH(t)\right)\left(1-\alpha_xy\int\frac{\underline{m}^2(z)t^2}{(1+t\underline{m}(z))^{2}}dH(t)\right)}dz\nonumber\\
    &&-\frac{\beta_x}{2\pi {\bf i}}\oint\limits_{\cal C}f_\ell(z)\frac{y\int\frac{\underline{m}^3(z)t^2}
      {(\underline{m}(z)t+1)^3}dH(t)}{1-y\int\frac{\underline{m}^2(z)t^2dH(t)}
      {(1+t\underline{m}(z))^2}}dz,\nonumber
  \end{eqnarray}
  and variance-covariance function
  \begin{eqnarray}
    &&{\rm Cov}(X_{f_{\ell'}},
    X_{f_\ell})\nonumber \label{covA3}\\
    &=&-\frac{1}{4\pi^2}\oint\limits_{{\cal
        C}_1}\oint\limits_{{\cal
        C}_2}\frac{f_{\ell'}(z_1)f_\ell(z_2)}{(\underline{m}(z_1)-\underline{m}(z_2))^2}
    \frac{\partial\underline{m}(z_1)}{\partial z_1}\frac{\partial\underline{m}(z_2)}{\partial z_2}dz_1dz_2\nonumber\\
    &&-\frac{y\beta_x}{4\pi^2}\oint\limits_{{\cal
        C}_1}\oint\limits_{{\cal
        C}_2}f_{\ell'}(z_1)f_{\ell}(z_2)\left[\int\frac{t}{(\underline{m}(z_1)t+1)^2}
      \frac{t}{(\underline{m}(z_2)t+1)^2}dH(t)\right] \frac{\partial\underline{m}(z_1)}{\partial z_1}\frac{\partial\underline{m}(z_2)}{\partial z_2}dz_1dz_2\nonumber\\
    &&{+}\frac{1}{4\pi^2}\oint\limits_{{\cal
        C}_1}\oint\limits_{{\cal
        C}_2} f_{\ell'}(z_1)f_{\ell}(z_2)\left[\frac{\partial^2}{\partial z_1\partial z_2}\log(1-a(z_1,z_2))\right]dz_1dz_2,\quad
  \end{eqnarray}
  where ${\cal C}$, ${\cal C}_1$ and ${\cal C}_2$ are closed contours in the complex plane enclosing  the support of
  the LSD $F^{y,H}$, and ${\cal C}_1$ and ${\cal C}_2$ are non-overlapping. Moreover, $a(z_1,z_2)$ is given by
  \[
  a(z_1,z_2)
  =
  \alpha_x\left(1+\frac{\underline{m}(z_1)\underline{m}(z_2)(z_1-z_2)}{\underline{m}(z_2)-\underline{m}(z_1)}\right).
  \]
\end{theorem}

\begin{remark}\label{remm}
  Theorem \ref{thm2}  gives   the explicit form for the   mean  and covariance  functions of $X_{f_\ell}$ expressed by contour integrals.
  For some specific  function $f_{\ell}$ and   $\bbQ$,
  the mean and covariance functions of  $X_{f_\ell}$ can be explicitly obtained.  Even if the mean and covariance functions  of  $X_{f_\ell}$  have no explicit forms, numerical methods can be used
  since the integrals are at most two-dimensional.  If the population
  is Gaussian, then $\beta_x=0$
  \mgai{and the corresponding terms in both the limiting mean and
    covariance functions vanish.}
  If the population is real, then we have  $\alpha_x=1$ and $$\frac{\partial^2\log(1-a(z_1,z_2))}{\partial z_1\partial z_2}=\frac{-1}{(\underline{m}(z_1)-\underline{m}(z_2))^2}
  \frac{\partial\underline{m}(z_1)}{\partial
    z_1}\frac{\partial\underline{m}(z_2)}{\partial z_2}.$$
  \mgai{Notice that  the limiting mean and covariance functions  of
    $\{X_{f_\ell}\}$
    depend on  the LSD $H$ of $\bbQ\bbQ^{*}$ and $y$ only.}
  In many
  applications, the ESD $H_n$ of $\bbQ\bbQ^{*}$ and the ratio
  $y_n=p/n$ can be used to replace the LSD $H$ and the limit $y$.
  But it is worthwhile to mention that when $\bbB_n$ is the centered sample covariance matrix $\bbB_n=(n-1)^{-1}\sum_{i=1}^{n}(\bby_i-\bar{\bby})(\bby_i-\bar{\bby})^T$
  with $\bar{\bby}=n^{-1}\sum_{i=1}^{n}\bby_i$, $y_n$ in (\ref{Yp}) should be replaced by $y_{n-1}$. For example, $\bby_i\sim N({\bf 0}_p, \bbI_p)$, then
  the Wishart matrices $n^{-1}\sum_{i=1}^{n}\bby_i\bby_i^T$ and $(n-1)^{-1}\sum_{i=1}^{n}(\bby_i-\bar{\bby})(\bby_i-\bar{\bby})^T$  have degrees of freedoms
  $n$ and $n-1$, respectively.
\end{remark}

For  practical use,
by Theorem \ref{thm2} and Remark \ref{rem001}, we will give three corollaries that are useful for computing the variance and covariance of the limiting distribution.
\begin{corollary}\label{coro0}
  Under the conditions of Theorem \ref{thm2} and if the population is real, then
  the random vector $(X_p(f_1),\ldots,X_p(f_L))$ converges to an $L$-dimensional Gaussian
  random vector $(X_{f_1},\ldots,X_{f_L})$ with mean function
  \begin{eqnarray*}
    \rE X_{f_\ell}&=&-\frac{1}{2\pi {\bf i}}\oint\limits_{\cal C}f_\ell(z)\frac{y\int\frac{\underline{m}^3(z)t^2}{(1+t\underline{m}(z))^{3}}dH(t)}
    {\left(1-y\int\frac{\underline{m}^2(z)t^2}{(1+t\underline{m}(z))^{2}}dH(t)\right)^2}dz-\frac{\beta_x}{2\pi {\bf i}}\oint\limits_{\cal C}f_\ell(z)\frac{y\int\frac{\underline{m}^3(z)t^2}{(\underline{m}(z)t+1)^3}dH(t)}{1-y\int\frac{\underline{m}^2(z)t^2dH(t)}{(1+t\underline{m}(z))^2}}dz,
  \end{eqnarray*}
  and variance-covariance function
  \begin{eqnarray*}
    &&{\rm Cov}(X_{f_{\ell '}}, X_{f_\ell})\\
    &=&-\frac{1}{2\pi^2}\oint\limits_{{\cal C}_1}\oint\limits_{{\cal
        C}_2}\frac{f_{\ell '}(z_1)f_\ell(z_2)}{(\underline{m}(z_1)-\underline{m}(z_2))^2} \frac{\partial\underline{m}(z_1)}{\partial z_1}\frac{\partial\underline{m}(z_2)}{\partial z_2}dz_1dz_2\\
    &&-\frac{y\beta_x}{4\pi^2}\oint\limits_{{\cal C}_1}\oint\limits_{{\cal C}_2}f_{\ell '}(z_1)f_{\ell}(z_2)\left[\int\frac{t}{(\underline{m}(z_1)t+1)^2}
      \frac{t}{(\underline{m}(z_2)t+1)^2}dH(t)\right] \frac{\partial\underline{m}(z_1)}{\partial z_1}\frac{\partial\underline{m}(z_2)}{\partial z_2}dz_1dz_2,
  \end{eqnarray*}
  where ${\cal C}$, ${\cal C}_1$ and ${\cal C}_2$ are
  closed contours in the complex plane enclosing  the support of
  the LSD $F^{y,H}$, and ${\cal C}_1$ and ${\cal C}_2$ are non-overlapping.
\end{corollary}

\begin{corollary}\label{coro1}
  Under the assumptions of Theorem \ref{thm2}, if the population is real and \mgai{$H_n=\delta_1$}, then we have
  $$
  F^{y,H_n}(f_\ell)=\int_{(1-\sqrt{y})^2}^{(1+\sqrt{y})^2}\frac{f_{\ell}(x)}{2\pi yx}\sqrt{[(1+\sqrt{y})^2-x][x-(1-\sqrt{y})^2]}dx+(1-1/y){\rm sgn}(y>1)f_{\ell}(0),
  $$
  \begin{eqnarray*}
    \rE X_{f_\ell}&=&\frac{f_{\ell}\left((1-\sqrt{y})^2\right)+f_{\ell}\left((1+\sqrt{y})^2\right)}{4} -
    \frac{1}{2\pi}\int_{(1-\sqrt{y})^2}^{(1+\sqrt{y})^2}\frac{f_{\ell}(x)}{\sqrt{4y-(x-1-y)^2}}dx\\
    &&-\frac{\beta_x}{2\pi {\bf i}}\oint\limits_{\cal C}f_\ell(z)\frac{y\underline{m}^3(z)
      (\underline{m}(z)+1)^{-3}}{1-y\underline{m}^2(z)(1+\underline{m}(z))^{-2}}dz,\nonumber
  \end{eqnarray*}
  and
  \begin{eqnarray*}
    {\rm Cov}(X_{f_{\ell '}},
    X_{f_\ell})\nonumber
    &=&-\frac{1}{2\pi^2}\oint\limits_{{\cal
        C}_1}\oint\limits_{{\cal
        C}_2}\frac{f_{\ell '}(z_1)f_\ell(z_2)}{(\underline{m}(z_1)-\underline{m}(z_2))^2}
    \frac{\partial\underline{m}(z_1)}{\partial z_1}\frac{\partial\underline{m}(z_2)}{\partial z_2}dz_1dz_2\\
    &&-\frac{y\beta_x}{4\pi^2}\oint\limits_{{\cal
        C}}\frac{f_{\ell'}(z_1)}{(\underline{m}(z_1)+1)^2}\frac{\partial\underline{m}(z_1)}{\partial z_1}dz_1\oint\limits_{{\cal
        C}}\frac{f_{\ell}(z_2)}{(\underline{m}(z_2)+1)^2}\frac{\partial\underline{m}(z_2)}{\partial z_2}dz_2,
  \end{eqnarray*}
  where ${\cal C}$, ${\cal C}_1$ and ${\cal C}_2$ are
  closed contours in the complex plane enclosing the support $[(1-\sqrt{y})^2, (1+\sqrt{y})^2]$ of
  the LSD $F^{y,H}$, ${\cal C}_1$ and ${\cal C}_2$ are non-overlapping, and $z=-\um^{-1}(z)+y(1+\um(z))^{-1}.$
\end{corollary}

If the real samples $\{\bby_i, i=1,\ldots,n\}$ are transformed as $\{(\bbQ\bbQ^{*})^{-1/2}\bby_i, i=1,\ldots,n\}$, then
the population covariance matrix of $(\bbQ\bbQ^{*})^{-1/2}\bby_i$ is
the identity matrix $\bbI_p$. Then \mgai{$H_n=\delta_1$}.
By Corollary \ref{coro1} and Remark \ref{remm}, we  give
\mgai{below a corollary useful for}  high-dimensional statistical inference.

\begin{corollary}\label{cor1}
  Under the assumptions of Theorem \ref{thm2}, if the population is
  such that  $H_n=\delta_1$, then  we have
  $$
  \left\{{\rm tr}[(\bbQ\bbQ^{*})^{-1}\bbB_n]-p,\ldots,{\rm tr}[(\bbQ\bbQ^{*})^{-1}\bbB_n]^L-pF^{(L)}(y_{n})\right\}^T\rightarrow N(\bm{\mu}, \bm{\Sigma})
  $$
  where $\bbB_n=n^{-1}\sum_{i=1}^n(\bby_i-\rE\bby_i)(\bby_i-\rE\bby_i)^{T},$
  $\bm{\mu}=(\mu_{1},\ldots,\mu_{L})^T$ and $\bm{\Sigma}=(\sigma_{\ell\ell'})_{\ell,\ell'=1}^{L}$. Moreover, if $\bbB_n=(n-1)^{-1}\sum_{i=1}^n(\bby_i-\bar{\bby})(\bby_i-\bar{\bby})^{T}$ with the sample mean $\bar{\bby}=n^{-1}\sum_{i=1}^{n}\bby_i$, then we have
  $$
  \left\{{\rm tr}[(\bbQ\bbQ^{*})^{-1}\bbB_n]-p,\ldots,{\rm tr}[(\bbQ\bbQ^{*})^{-1}\bbB_n]^L-pF^{(L)}(y_{n-1})\right\}^T\rightarrow N(\bm{\mu}, \bm{\Sigma}).
  $$
  Here,
  \begin{eqnarray*}
    F^{(\ell)}(y)&=&\int_{(1-\sqrt{y})^2}^{(1+\sqrt{y})^2}\frac{x^{\ell-1}}{2\pi y}\sqrt{[(1+\sqrt{y})^2-x][x-(1-\sqrt{y})^2]}dx,~~\ell=1,\ldots,L,
  \end{eqnarray*}
  \begin{eqnarray*}
    \mu_{\ell}&=&\frac{(1-\sqrt{y})^{2\ell}+(1+\sqrt{y})^{2\ell}}{4}-\frac{1}{2}\sum\limits_{\ell_1=0}^{\ell}\binom{\ell}{\ell_1}^2y^\ell\\
    &&+\beta_x{\rm sgn}(\ell\geq 2)\sum\limits_{\ell_2=2}^\ell\binom{\ell}{\ell_2-2}\binom{\ell}{\ell_2}y^{\ell+1-\ell_2},~~\ell=1,\ldots,L,\end{eqnarray*}
  and
  \begin{eqnarray*}
    \lefteqn{\sigma_{\ell\ell'}=}\\
    &&2y^{\ell+\ell'}\sum\limits_{\ell_1=0}^{\ell-1}\sum\limits_{\ell_2=0}^{\ell'}
    \binom{\ell}{\ell_1}\binom{\ell'}{\ell_2}\left(\frac{1-y}{y}\right)^{\ell_1+\ell_2}
    \sum\limits_{\ell_3=0}^{\ell-\ell_1}\ell_3\binom{2\ell-1-\ell_1-\ell_3}{\ell-1}\binom{2\ell'-1-\ell_2+\ell_3}{\ell'-1}\\
    &&+y\beta_x\sum\limits_{\ell_3=1}^{\ell}\binom{\ell}{\ell_3-1}\binom{\ell}{\ell_3}y^{\ell-\ell_3}\sum\limits_{\ell_3=1}^{\ell'}\binom{\ell'}{\ell_3-1}
    \binom{\ell'}{\ell_3}y^{\ell'-\ell_3},~~\ell, \ell'=1,\ldots,L.
  \end{eqnarray*}\label{corollary1}
\end{corollary}
In practice,  $y$ may be replaced by $y_{n-1}$. Corollary
\ref{corollary1} can be used for many high dimensional testing problems
and its proof is given in the appendix.






\section{Applications}
\label{sec:app}
In this section, we
\mgai{propose two applications of
  Theorem~\ref{thm2}.  The first application is about testing the
  structure of a high-dimensional covariance matrix.
  The second application applies such structure testing to the
  analysis of a  large fMRI dataset.}

\subsection{Testing  high dimensional covariance structure}
\label{ssec:test}
Testing  the structure of
high dimensional covariance matrices is  \mgai{a fundamental problem}
in multivariate and high-dimensional data analysis
\mgai{\citep{Sriv05, bjyz09, CZZ2010, WY13}.
  Most of this literature consider samples of the form
  $\bby_i=\bbQ\bbx_i$ where the dimension $k$ of the  $\bbx_i$'s is
  finite.
  When $k=\infty$ as for time series observations, these
  existing results are not applicable anymore for structure  testing
  on   high dimensional covariance matrices.}

We consider two     testing problems on high dimensional covariance structure. Let $\bm{\Sigma}_0$ be  a given covariance matrix.
The first  testing problem is given by
\begin{equation}\label{H01}
  H_{01}: \bbQ\bbQ^{*}=\bm{\Sigma}_0.
\end{equation}
We consider
$\rtr(\bm{\Sigma}^{-1/2}_0\bbB_n\bm{\Sigma}_0^{-1/2}-\bbI_p)^2$ as our test statistic,
where
$\bbB_n=(n-1)^{-1}\sum_{i=1}^{n}(\bby_i-\bar{\bby})(\bby_i-\bar{\bby})^T$. Because
$\bm{\Sigma}_0^{-1/2}\bbQ\bbQ^{*}\bm{\Sigma}_0^{-1/2}=\bbI_p$ under
$H_{01}$,  \mgai{$H_n=\delta_1$ (Dirac mass at 1)}. Therefore, it
follows from Corollary \ref{cor1} and the delta method that  we have
under $H_{01}$
\begin{equation}
  \frac{1}{2}\{\rtr(\bm{\Sigma}^{-1}_0\bbB_n-\bbI_p)^2-py_{n-1}-(\beta_x+1)y_{n-1}\}/\sqrt{y^2_{n-1}+(\beta_x+2)y^3_{n-1}}
  ~\mgai{\Longrightarrow}~
  N(0, 1).
\end{equation}
The second testing \mgai{problem we consider is about the hypothesis}
\begin{equation}
  H_{02}:  \bbQ\bbQ^{*}=\sigma^2\bm{\Sigma}_0.
\end{equation}
\mgai{When $\bm{\Sigma}_0=\bbI_p$, this reduces to the well known
  sphericity test \citep{WY13}.
}
We consider
${\rm tr}\{\bm{\Sigma}_0^{-1/2}\bbB_n\bm{\Sigma}_0^{-1/2}/[p^{-1}{\rm tr}(\bm{\Sigma}_0^{-1}\bbB_n)]-\bbI_p\}^2$ as our test statistic.
Since $\sigma^{-2}\bm{\Sigma}_0^{-1/2}\bbQ$
$\bbQ^{*}\bm{\Sigma}_0^{-1/2}=\bbI_p$ under $H_{02}$,  \mgai{again $H_n=\delta_1$}.
By Corollary \ref{cor1} and the delta method, we have under $H_{02}$,
\begin{equation}\label{H02}
  \frac{1}{2}\{{\rm tr}[\bm{\Sigma}_0^{-1}\bbB_n/(p^{-1}{\rm
    tr}(\bm{\Sigma}_0^{-1}\bbB_n))-\bbI_p]^2-py_{n-1}-(\beta_x+1)y_{n-1}\}/y_{n-1}
  ~\mgai{\Longrightarrow}~ N(0, 1).
\end{equation}

Simulation experiments are done to evaluate the finite sample
performance of our test \mgai{for the null} $H_{02}$.
The test size is set at $5\%$ and
\mgai{empirical sizes and powers are obtained using 5000 independent replications}.
We draw $\bby_i$ according to the AR(2) model $y_{ti}=\phi_1y_{t-1,i}+\phi_2y_{t-2,i}+e_{ti}$,
where the $e_{ti}$'s are i.i.d. as  $N(0, \sigma^2)$.
We set $n\in \{100, 200,  300\}$ and $p\in \{50, 100, 200, 500, 1000\}$.
\mgai{Empirical sizes are evaluated for $\phi_1\in\{0.3, 0.6\}$, $\phi_2\in\{0.2, 0.3\}$}, and
$\bm{\Sigma}_0=\{\gamma_{|i-j|}\}_{i, j=1}^{p}$ with $\gamma_0=1$, $\gamma_1=\phi_1/(1-\phi_2)$ and
$\gamma_{\ell}=\phi_1\gamma_{\ell-1}+\phi_2\gamma_{\ell-2}$ for
$\ell\geq 2$.
\mgai{They are given} in  Table~\ref{table1}.
To evaluate empirical powers, the null hypothesis has $\phi_1=\phi_2=0.18$ and the alternative hypothesis 
has $\phi_1\in\{0.3, 0.35\}$ and $\phi_2\in\{0.2, 0.25\}$. Table~\ref{table2} presents these  empirical powers.
\mgai{Both tables show a very satisfactory  finite-sample performance
  of the proposed test.}

{\small
\begin{table}[!h]
  \caption{Empirical test sizes for $H_{02}$ (in percentage)}
  \doublerulesep 0.5pt
  \begin{center}
    \begin{tabular}{ccr|ccccc|ccccc}
      \hline
      \hline
      & & & \multicolumn{5}{c|}{$\phi_1=0.3$} & \multicolumn{5}{c}{$\phi_1=0.6$}\cr
      $\phi_2$ & $n$  & & $p=$50 & 100 & 200 & 500 & 1000 & p=50 & 100 & 200 & 500 & 1000\cr
      \hline
      0.2          & 100 &&5.36&5.38&5.38&5.32&5.54&5.12&5.36&5.92&5.46&4.94\cr
      & 200 &&5.12&5.46&5.22&5.22&5.80&5.76&5.34&5.32&5.14&4.76\cr
      & 300 &&5.12&4.82&5.10&5.50&4.90&5.44&5.22&5.36&5.52&5.10 \cr
      \hline
      0.3          & 100 &&5.72&5.28&5.60&5.22&5.90&4.90&5.76&5.42&5.06&5.28 \cr
      & 200 &&4.60&5.12&5.28&5.02&5.16&5.52&5.68&5.48&5.18&5.20 \cr
      & 300 &&5.04&5.56&5.30&5.24&4.88&5.54&5.64&5.48&5.50&5.56 \cr
      \hline
    \end{tabular}
  \end{center}\label{table1}
\end{table}
}

{\small
\begin{table}[!h]
  \caption{Empirical powers for $H_{02}$ (in percentage)}
  \doublerulesep 0.5pt
  \begin{center}
    \begin{tabular}{ccr|ccccc|ccccc}
      \hline
      \hline
      & & & \multicolumn{5}{c|}{$\phi_1=0.3$} & \multicolumn{5}{c}{$\phi_1=0.35$}\cr
      $\phi_2$ & $n$  & & $p=$50 & 100 & 200 & 500 &  & p=50 & 100 & 200 & 500 &\cr
      \hline
      0.2          & 100 &&59.00&61.24&62.80&64.50&&98.54&99.36&99.54&100.00\cr
      & 200 &&97.98&98.90&99.28&100.00&&100.00&100.00&100.00&100.00\cr
      & 300 &&100.00&100.00&100.00&100.00&&100.00&100.00&100.00&100.00 \cr
      \hline
      0.25     & 100&&94.46&97.18&97.46&98.3&&100.00&100.00&100.00&100.00\cr
      & 200&&100.00&100.00&100.00&100.00&&100.00&100.00&100.00&100.00\cr
      & 300&&100.00&100.00&100.00&100.00&&100.00&100.00&100.00&100.00\cr
      \hline
    \end{tabular}
  \end{center}\label{table2}
\end{table}
}

\subsection{Application to a fMRI dataset}
\label{sec:fMRI}

We apply the proposed method  to \mgai{a structure testing problem for}
the resting-state fMRI data from the ADHD-200 sample. The whole data
was collected from eight sites of the ADHD-200 consortium. We focus
on the data from New York University (NYU) with the largest number
of \mgai{subjects. Among them,} we used the 1000 ROI extracted time courses preprocessed
by the Neuroimaging Analysis Kit (NIAK) \citep{Ahn2015}. The data set
consists of   70 subjects $(n=70)$, each of which contains 924 Regions
Of Interest (ROIs). Moreover, each ROI contains a time series with 172
time points. At a given ROI, the fMRI observations   for the
$i$th subject are treated as $\bby_i$,  \mgai{a time series of length $p=172$.}



\begin{figure}[htbp]
  \begin{center}
    \includegraphics[height=3in,width=3in, angle=0]{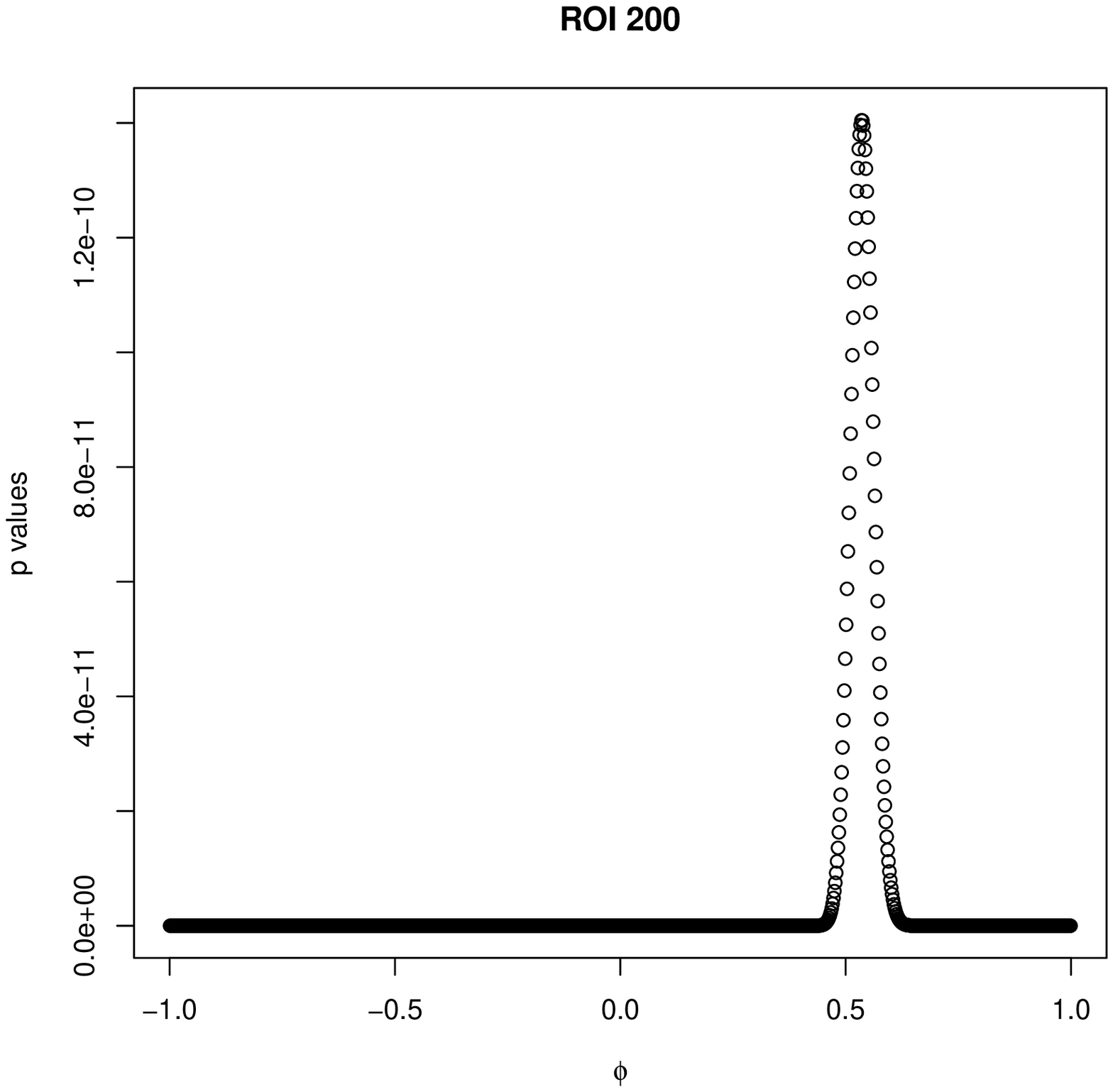}\quad
    \includegraphics[height=3in,width=3in, angle=0]{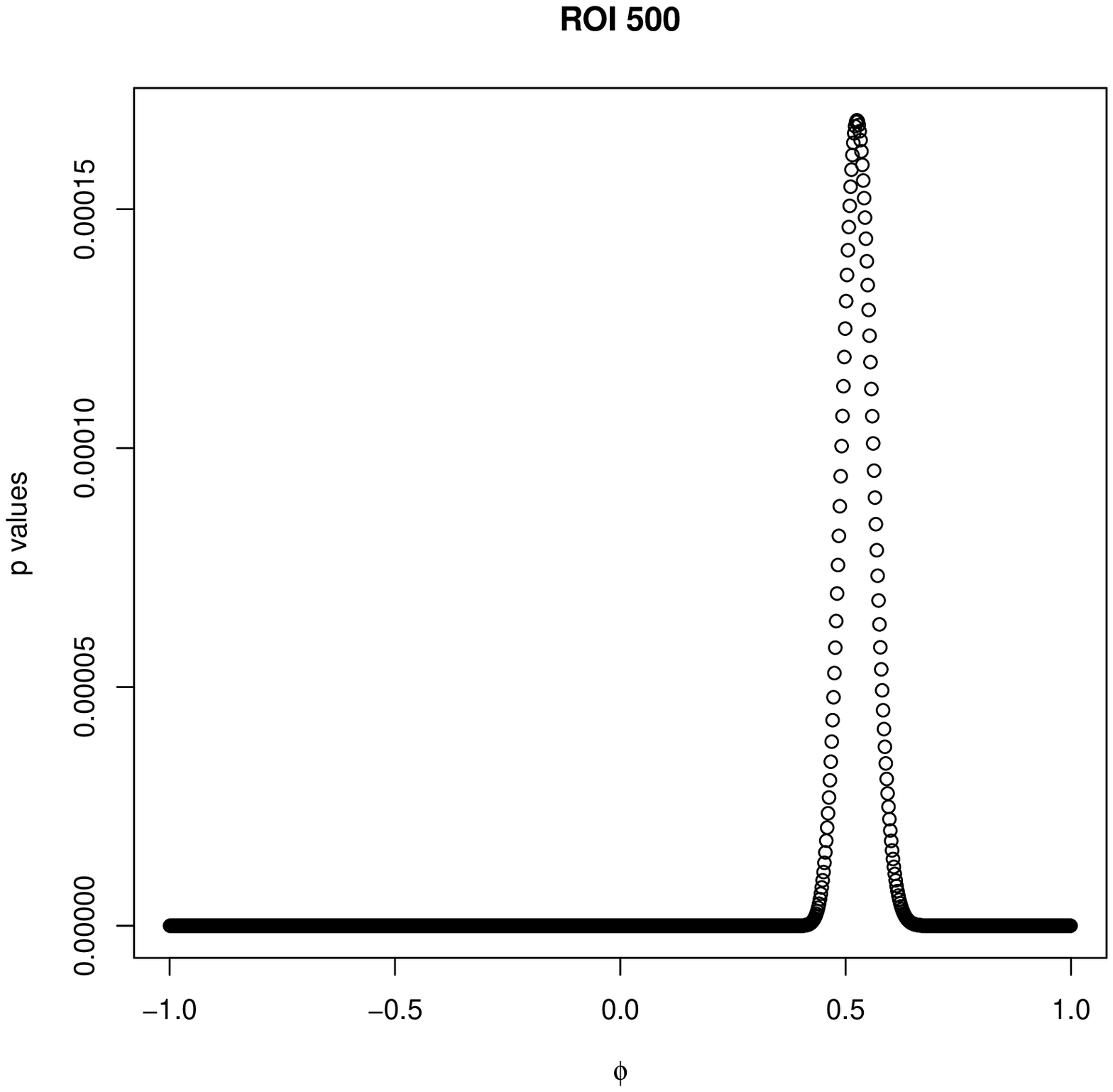}
  \end{center}
  \vskip-1cm
  \caption{P-values of both ROI 200 and 500 for the testing problem (\protect\ref{test1}).\label{fig2}}
\end{figure}

The question of interest here is to test whether the standard
autoregressive assumption (e.g., AR(1) or AR(2)) is valid for the
processed fMRI time series across all ROIs \citep{Lindquist}.
\mgai{Applying results derived in the previous section, we first test the AR(1) structure with the hypothesis}
\begin{equation}\label{test1}
  H_0: \mbox{The time series in a specific ROI is an AR(1) process}.
\end{equation}
Let $\phi$ be the coefficient of the AR(1) process and $\bm{\Sigma}_{\phi}=(\phi^{|i-j|})_{i, j=1}^{p}$. The auto-covariance matrix of the AR(1) process is  given by $\bbQ\bbQ^{*}=\sigma^2\bm{\Sigma}_{\phi}$ and then we have ${\rm Cov}(\bm{\Sigma}_{\phi}^{-1/2}{\bf y}_i)=\sigma^2{\bf I}_p$, where  $\sigma^2$ is the variance.
Testing the AR(1) assumption of a given ROI is equivalent to testing
the sphericity covariance structure of such ROI. Similar to the
arguments of Section~\ref{ssec:test},  we can show \mgai{that
  (\ref{H02}) still holds here}.  \mgai{We detail results for ROI 200  and 500}.
Since $\phi$ takes values in the interval $(-1, 1)$,  we set $\phi$ be
$-1+0.001, -1+0.002,\ldots,1-0.001$ and calculated the P-values of the
test (\ref{test1}) in Figure~\ref{fig2}. It indicates that both
\mgai{sets of} P-values are much smaller than $5\%$, which yields the
rejection of the proposed  AR(1) structure in both ROIs.

\mgai{Next we test the AR(2) hypothesis}
\begin{equation}\label{test2}
  H_0: \mbox{The time series in a given ROI   is an  AR(2) process}.
\end{equation}
Let $\phi_1$ and $\phi_2$ be the auto-regressive coefficients of the
AR(2) model. Then, the auto-covariance matrix of the AR(2) model is
equal to $\bbQ\bbQ^{*}=\sigma^2\bm{\Sigma}_{\phi_1, \phi_2}$, where
$\bm{\Sigma}_{\phi_1, \phi_2}=(\gamma_{|i-j|})_{i, j=1}^{p}$ with
$\gamma_0=1$, $\gamma_1=\phi_1/(1-\phi_2)$, and
$\gamma_{\ell}=\phi_1\gamma_{\ell-1}+\phi_2\gamma_{\ell-2}$ for
$\ell\geq 2$ satisfying $\phi_1^2+\phi_2^2<1$ and
$\phi_2+|\phi_1|<1$. Then,  we have ${\rm Cov}(\bm{\Sigma}_{\phi_1,
  \phi_2}^{-1/2}{\bf y}_i)=\sigma^2{\bf I}_p$. For a specific ROI,
testing  the AR(2) assumption is again equivalent to  testing the
sphericity covariance structure;
\mgai{in particular (\ref{H02}) still holds here.}
Similar to the previous test,
we calculated the P-values of our test statistic for every $\phi_1,
\phi_2\in\{-1+0.001, -1+0.002, \ldots, 1-0.001\}$ satisfying
$\phi_1^2+\phi_2^2<1$ and $\phi_2+|\phi_1|<1$ for ROI 200 and 500.
\mgai{As a result},  these  P-values from both ROI 200 and 500 are
much smaller than the nominal level $5\%$ so that  the proposed
$AR(2)$ structure is clearly rejected.
\mgai{In summary,
  neither of the AR(1) and AR(2) structures is}
suitable for the resting-state fMRI data considered here.

\section{Proofs of the main theorems}\label{sec:proofs}

\subsection{Proof of Theorem \protect\ref{thm1}}

\subsubsection{\mgai{Case where  $\bbQ$ has an infinite number of  columns}}\label{511}

{Consider  a sequence of $\{k_p\}$ that satisfies $p^{-1}\sum_{i=1}^p\sum_{\ell=k_p+1}^{\infty}|q_{i\ell}^2|=o(p^{-2})$.
  For simplicity, we write $k$ for $k_p$.
  If $\bbQ$ is a $p\times\infty$ dimensional matrix, then we truncate
  $\bbQ$ as
  $\bbQ=(\hat{{\bbQ}}, \tilde{{\bbQ}})$
  where $\hat{{\bbQ}}$ is a $p\times k$ dimensional matrix and $\tilde{{\bbQ}}=(q_{ij})$
  is a $p\times \infty$ dimensional matrix with $i=1,\ldots,p, j=k+1,\ldots,\infty$.
  Similarly, truncate $\bbX_n$ as $\bbX_n=(\hat{{\bbX}}_n, \tilde{{\bbX}}_n)$
  where $\hat{{\bbX}}_n$ is $k\times n$ dimensional and $\tilde{{\bbX}}_n$ is $\infty\times n$ dimensional.
  Then we have
  $$
  L^4(F^{n^{-1}\bbQ\bbX_n\bbX_n^*\bbQ^*}, F^{n^{-1}\hat{{\bbQ}}\hat{{\bbX}}_n\hat{{\bbX}}_n^*\hat{{\bbQ}}^*})
  \leq2p^{-2}n^{-2}{\rm tr}(\bbQ\bbX_n\bbX_n^*\bbQ^*+\hat{{\bbQ}}\hat{{\bbX}}_n\hat{{\bbX}}_n^*\hat{{\bbQ}}^*)
  {\rm tr}\tilde{{\bbQ}}\tilde{{\bbX}}_n\tilde{{\bbX}}_n^*\tilde{{\bbQ}}^*
  $$
  where $L(, )$ is the Levy distance, $F^{n^{-1}\bbQ\bbX_n\bbX_n^*\bbQ^*}$ is the ESD of $n^{-1}\bbQ\bbX_n\bbX_n^*\bbQ^*$ and
  $F^{n^{-1}\hat{{\bbQ}}\hat{{\bbX}}_n\hat{{\bbX}}_n^*\hat{{\bbQ}}^*}$ is the ESD of $n^{-1}\hat{{\bbQ}}\hat{{\bbX}}_n\hat{{\bbX}}_n^*\hat{{\bbQ}}^*$.
  We have
  \begin{eqnarray*}
    &&\frac1{pn}\rtr(\tilde{{\bbQ}}\tilde{{\bbX}}_n\tilde{{\bbX}}_n^*\tilde{{\bbQ}}^*)=\frac1{pn}\sum_{i=1}^p\sum_{j=1}^n\left|\sum_{\ell=k+1}^{\infty}
      q_{i\ell}x_{\ell j}\right|^2\non
    &=&\frac1{pn}\sum_{i=1}^p\sum_{j=1}^n\sum_{\ell=k+1}^{\infty}|q_{i\ell}^2||x_{\ell j}^2|+\frac1{pn}\sum_{i=1}^p\sum_{j=1}^n\sum_{k_1\ne k_2}q_{ik_1}\bar q_{ik_2} x_{k_1,j}\bar{x}_{k_2,j}.
  \end{eqnarray*}
  Note that
  $$
  \rE\left(\frac1{pn}\sum_{i=1}^p\sum_{j=1}^n\sum_{\ell=k+1}^{\infty}|q_{ik}^2||x_{\ell j}^2|\right)=\frac1{p}\sum_{i=1}^p\sum_{\ell=k+1}^{\infty}|q_{i\ell}^2|=o(p^{-2})
  $$
  and
  \bqa
  &&\rE\left(\frac1{pn}\sum_{i=1}^p\sum_{j=1}^n\sum_{\ell=k+1}^{\infty}|q_{i\ell}^2|(|x_{\ell j}^2|-\rE|x_{\ell j}^2|)\right)^4\non
  &\le&\frac1{p^4n^4}\sum_{j=1}^n\sum_{\ell=k+1}^{\infty}\left(\sum_{i=1}^p|q_{i\ell}^2|\right)^4\rE|x_{\ell j}^8|
  +\frac3{p^4n^4}\left(\sum_{j=1}^n\sum_{\ell=k+1}^{\infty}\left(\sum_{i=1}^p|q_{i\ell}^2|\right)^2\rE|x_{\ell j}^4|\right)^2 \non
  &\le&\frac{\eta_n^6}{p^4}\sum_{\ell=k+1}^{\infty}\|\bbq_\ell\|^2+\frac{3\eta_n^4}{p^4}\left(\sum_{\ell=k+1}^{\infty}\|\bbq_\ell\|^2\right)^2=o(p^{-2}).
  \eqa
  These inequalities simply imply $(pn)^{-1}\sum_{i=1}^p\sum_{j=1}^n\sum_{\ell=1}^k|q_{i\ell}^2||x_{\ell j}^2|=o_{a.s.}(1)$.
  Furthermore,
  \bqa
  &&\rE\left(\frac1{pn}\sum_{i=1}^p\sum_{j=1}^n\sum_{k_1\ne k_2}q_{ik_1}\bar q_{ik_2}x_{k_1,j}\bar{x}_{k_2,j}\right)^2
  \le \frac2{p^2n^2}\sum_{j=1}^n\sum_{k_1\ne k_2}\left|\sum_{i=1}^pq_{ik_1}\bar q_{ik_2}\right|^2\non
  &\le&\frac2{p^2n}\rtr (\tilde{{\bbQ}}\tilde{{\bbQ}}^*)^2=o(p^{-2}),
  \eqa
  which implies that $(pn)^{-1}\sum_{i=1}^p\sum_{j=1}^n\sum_{k_1\ne k_2}q_{ik_1}\bar q_{ik_2}x_{k_1,j}\bar{x}_{k_2,j}\to0,\,a.s.$ Then we have
  $$(pn)^{-1}\rtr(\tilde{{\bbQ}}\tilde{{\bbX}}_n\tilde{{\bbX}}_n^*\tilde{{\bbQ}}^*)=o_{a.s.}(1).$$
  Similarly, we can prove $(pn)^{-1}{\rm tr}(\bbQ\bbX_n\bbX_n^*\bbQ^*)=O_{a.s.}(1)$ and  $(pn)^{-1}{\rm tr}(\hat{{\bbQ}}\hat{{\bbX}}_n\hat{{\bbX}}_n^*\hat{{\bbQ}}^*)=O_{a.s.}(1)$.
  Then we have
  $$
  L^4(F^{n^{-1}\bbQ\bbX_n\bbX_n^*\bbQ^*}, F^{n^{-1}\hat{{\bbQ}}\hat{{\bbX}}_n\hat{{\bbX}}_n^*\hat{{\bbQ}}^*})=o_{a.s.}(1).
  $$
  Therefore, without loss of generality, we will
  \mgai{hereafter
    assume that the number of columns $k$ of
    $\bbQ$ is finite.}}

\subsubsection{\mgai{Sketch}  of the proof of Theorem \ref{thm1}}
\cite{Silv95} obtained the LSD of the classical sample covariance matrix when the dimension increases proportionally with the sample size.
But the sample covariance matrix ${\bf B}_n$ of this paper is different from that of \cite{Silv95} because \cite{Silv95} assumed that ${\bf Q}$ is $p\times p$
and nonnegative definite while this paper  assumes ${\bf Q}$ is $p\times k$ with $k\ge p$, even being infinity. In establishing the LSD of ${\bbB_n}$ under the assumptions of Theorem \ref{thm1}, we need the truncation, centralization and rescaling on $x_{ij}$. Under the assumptions of Theorem \ref{thm1}, the LSD of ${\bf B}_n$ will be the same before or after the truncation, centralization and rescaling. Because the dimension of ${\bf Q}$ is different from those defined in \cite{Silv95}, the assumptions and truncations of Theorem \ref{thm1} are also different from those given in \cite{Silv95}. The readers are reminded that the corresponding proofs given in \cite{BS10} and \cite{Silv95}, strongly depend on a fact that $\sum_{j=1}^n\beta_j(z)=\rtr(n^{-1}\bbX_n^*\tilde \bbT_n\bbX_n-z\bbI_n)^{-1}$, where $\beta_j(z)$'s are defined in Subsection \ref{subsec511} and $\tilde\bbT_n=\bbT_n$ in \cite{Silv95} while $\tilde\bbT_n=\bbQ^*\bbQ$ in the present paper. It may cause some ambiguity when $k$ is infinity in the multiplication of infinite-dimensional matrices, we will avoid to use this fact in this paper.

Let $m_n(z)=p^{-1}\rtr(\bbB_n-z\bbI_p)^{-1}$ be the Stieltjes transform of the ESD of $\bbB_n$, where $z=u+iv$ with $v>0$. We
will show that $\um_n(z)$ tends to a non-random limit that satisfies (\ref{eq1}), where $\um_n(z)=-(1-y_n)z^{-1}+y_n m_n(z)$ is in fact the Stieltjes transform of $\underline \bbB_n=n^{-1}\bbX_n^*\widetilde{\bbT}_n\bbX_n$, where $\widetilde{\bbT}_n=\bbQ^*\bbQ$. The proof will be split into
two parts: \bqa
\mbox{(i)}&&m_n(z)-\rE m_n(z)\to 0~a.s.,\label{eqg1}\\
\mbox{(ii)}&&\rE\um_n(z)\to\um(z)~\mbox{satisfying
  (\ref{eq1})}.\label{eqg2} \eqa For truncation and re-normalization (see Appendix B), we may further assume $|x_{ij}|<\eta_n\sqrt{n}/\|\bbq_i\|$, where the constant sequence $\{\eta_n\}$ tends to zero as $n\to\infty$.

\subsubsection{Proof of (\ref{eqg1})}\label{subsec511}

Let $\bbr_j=n^{-1/2}\bbQ\bbx_j$, then $\bbB_{n}=\sjln\bbr_j\bbr_j^*$.
Define \mgai{$ \bbD(z)=\bbB_n-z\bbI_p$ and }
\begin{equation}
  \label{eq:betaj-gammaj}
  \bbD_{j}(z)=\bbD(z)-\bbr_j\bbr_j^*,\quad
  \beta_j(z)=[1+\bbr_j^*\bbD_{j}^{-1}(z)\bbr_j]^{-1},\quad
  \bgma_j(z)=\bbr_j^*\bbD_j^{-2}(z)\bbr_j\beta_j(z).
\end{equation}
Since $\Im(\beta_j^{-1}(z))=v\bbr_j^*\bbD_{j}^{-1}(z)(\bbD_{j}^{-1}(z))^*\bbr_j>v\Big|\bbr_j^*\bbD_{j}^{-2}(z)\bbr_j\Big|$,
we obtain \be |\bgma_j(z)|\le v^{-1},\label{eqgp1}\ee where $\Im(z)=v$.
Moreover,
\begin{eqnarray*}
  \Im [z\bbr_j^{*}\bbD_j^{-1}(z)\bbr_j]
  &=&(2{\bf i})^{-1}[z\bbr_j^{*}\bbD_j^{-1}(z)\bbr_j-\bar{z}\bbr_j^{*}\overline{\bbD_j^{-1}}(z)\bbr_j]\\
  &=&v|z|^{-2}z\bbr_j^{*}\bbD_j^{-1}(z)\left(\sum\limits_{i\not=j}\bar{z}\bbr_i\bbr_i^{*}\right)\overline{\bbD_j^{-1}}(z)\bbr_j\geq 0,
\end{eqnarray*}
where $\bar{z}$ and $\overline{\bbD_j^{-1}}(z)$ denote the conjugate of $z$ and $\bbD_j^{-1}(z)$. Thus we have
\begin{equation}\label{betaj}
  |\beta_j(z)|\leq |z|v^{-1}.
\end{equation}
Denote the conditional expectation given $\{\bbr_1,\cdots,\bbr_j\}$ by $\rE_j$ and $\rE_0$ for the
unconditional expectation. Then, we have \bqa m_n(z)-\rE m_n(z)
&=&p^{-1}\sjln (\rE_j-\rE_{j-1})\rtr\bbD^{-1}(z)\non
&=&p^{-1}\sjln (\rE_j-\rE_{j-1})[\rtr\bbD^{-1}(z)-\rtr\bbD_{j}^{-1}(z)]\non
&=&p^{-1}\sjln (\rE_j-\rE_{j-1})\bgma_j(z).\label{eqgp2} \eqa By
Lemma \ref{lem1} (Burkholder inequality) and the inequality (\ref{eqgp1}), for any $\ell>1$,
we obtain \bqa\rE|m_n(z)-\rE m_n(z)|^{\ell}&\le&
K_0p^{-\ell}\rE\left(\sjln |(\rE_j-\rE_{j-1})\bgma_j(z)|^2\right)^{\frac{\ell}{2}}\le K_0v^{-\ell}p^{-\ell}n^{\frac{\ell}{2}}\label{eqgp3}\qquad\eqa where $K_0$ is a constant. Taking $\ell>2$, Chebyshev inequality and (\ref{eqgp3}) imply (\ref{eqg1}): $m_n(z)-\rE m_n(z)\to 0$ a.s.

\subsubsection{Proof of (\ref{eqg2})}
Following the steps of the proof of Theorem 1.1 of \cite{BZ2008}, define $\bbK=(1+y_na_{n,1})^{-1}\bbT_n$ and $y_n=p/n$,
where $\bbT_n=\bbQ\bbQ^*$, $a_{n,\ell}=p^{-1}\rE\rtr[\bbT_n^{\ell}\bbD^{-1}(z)]$, $\ell=0$ or
$1$. We have
\bqn
(\bbK-z\bbI)^{-1}-\bbD^{-1}(z)&=&\sum\limits_{j=1}^n(\bbK-z\bbI)^{-1}\bbr_{j}\bbr_{j}^*\bbD^{-1}(z)-(\bbK-z\bbI)^{-1}\bbK\bbD^{-1}(z)\\
&=&\sum\limits_{j=1}^n(\bbK-z\bbI)^{-1}\bbr_{j}\bbr_{j}^*\bbD^{-1}_{j}(z)\beta_{j}(z)-(\bbK-z\bbI)^{-1}\bbK\bbD^{-1}(z).
\eqn
For $\ell=0,1$, multiplying both sides by $\bbT_n^{\ell}$ and then taking trace and dividing by $p$, we have \bqa
&&p^{-1}\rE\rtr[\bbT_n^{\ell}(\bbK-z\bbI)^{-1}]-a_{n,\ell}\non
&=&p^{-1}\sum_{j=1}^n\rE\bbr_j^*\bbD^{-1}_j(z)\bbT_n^{\ell}(\bbK-z\bbI)^{-1}\bbr_j\beta_j(z)-p^{-1}\rE\rtr
[\bbT_n^{\ell}(\bbK-z\bbI)^{-1}\bbK\bbD^{-1}(z)]. \label{eq0011} \eqa
One can prove a formula similar to  (1.15) of \cite{BS04} and verify that
\bqa
&&\rE\left|\bbr_j^*\bbD_j^{-1}(z)\bbr_j-n^{-1}\rtr (\bbQ^*\bbD_j^{-1}(z)\bbQ)\right|^2\non
&\le&2n^{-2}\rE\rtr(\bbQ^*\bbD_j^{-1}(z)\bbQ)(\bbQ^*\bbD_j^{-1}(z)\bbQ)^*
+n^{-2}\sum_{i=1}^k |(\bbQ^*\bbD_j^{-1}(z)\bbQ)_{ii}|^2\rE|x_{ij}^4|\non
&\le&2n^{-2}pv^{-2}\|\bbQ\|^4+n^{-2}v^{-2}\sum_{i=1}^k\|\bbq_i\|^4\eta_n^2n/
\|\bbq_i\|^2\le C\eta_n^2\to0.
\label{eq3501}
\eqa
By Lemma \ref{lem7}, one can prove that
$$
\left|n^{-1}\rtr (\bbQ^*\bbD_j^{-1}(z)\bbQ)-n^{-1}\rtr (\bbQ^*\bbD^{-1}(z)\bbQ)\right|^2\le K(n^2v^2)^{-1}
$$
and by the similar method of the proof of (\ref{eqg1}), we have
$$
\rE\left|n^{-1}\rtr (\bbQ^*\bbD^{-1}(z)\bbQ)-n^{-1}\rE\rtr (\bbQ^*\bbD^{-1}(z)\bbQ)\right|^2\le K(nv^2)^{-1}.
$$
Therefore we have
$\rE|\beta_j^{-1}(z)-(1+y_n a_{n,1})|^2=o(1)
\label{eq350}$
which, applying Lemma \ref{lem7} again and (\ref{betaj}), implies that
\bqa
&&p^{-1}\sum_{j=1}^n\rE\bbr_j^*\bbD^{-1}_j(z)\bbT_n^{\ell}(\bbK-z\bbI)^{-1}\bbr_j\beta_j(z)-p^{-1}\rE\rtr[\bbT_n^{\ell}(\bbK-z\bbI)^{-1}\bbK\bbD^{-1}(z)]\non
&=&p^{-1}\sum_{j=1}^n\rE\bbr_j^*\bbD^{-1}_j(z)\bbT_n^{\ell}(\bbK-z\bbI)^{-1}\bbr_j(1+y_na_{n,1})^{-1}\non
&&-p^{-1}\rE\rtr[\bbT_n^{\ell}(\bbK-z\bbI)^{-1}\bbK\bbD^{-1}(z)]+o(1)\non
&=&(pn)^{-1}\sum_{j=1}^n\rE\rtr\bbD^{-1}_j(z)\bbT_n^{\ell}(\bbK-z\bbI)^{-1}\bbT_n(1+y_na_{n,1})^{-1}\non
&&-p^{-1}\rE\rtr[\bbT_n^{\ell}(\bbK-z\bbI)^{-1}\bbK\bbD^{-1}(z)]+o(1)\non
&=&(pn)^{-1}\sum_{j=1}^n\rE\rtr\bbD^{-1}_j(z)\bbT_n^{\ell}(\bbK-z\bbI)^{-1}\bbK-p^{-1}\rE\rtr \bbD^{-1}(z)\bbT_n^{\ell}(\bbK-z\bbI)^{-1}\bbK+o(1)\non&=&o(1).\label{eq0012}
\eqa
It then follows from (\ref{eq0011}) and (\ref{eq0012}) that
\bqa
a_{n,\ell}&=&p^{-1}\rtr\Big\{\bbT_n^{\ell}\big[(1+y_na_{n,1})^{-1}\bbT_n-z\bbI\big]^{-1}\Big\}+o(1)\non
&=&\int\frac{t^\ell}{t(1+y_na_{n,1})^{-1}-z}dH_n(t)+o(1),
\label{eq0013}
\eqa
where $H_n$ is the ESD of $\bbT_n$. Because $\Im[z(1+y_na_{n,1})]>v$, we conclude that $|(1+y_na_{n,1})^{-1}|\le |z|/v$. Taking $\ell=1$ in (\ref{eq0013}) and multiplying both sides by $(1+y_na_{n,1})^{-1}$, we obtain
\bqn
\frac{a_{n,1}}{1+y_na_{n,1}}&=&\int\frac{t(1+y_na_{n,1})^{-1}}{t(1+y_na_{n,1})^{-1}-z}dH_n(t)+o(1)\non
&=&1+z\int\frac{1}{t(1+y_na_{n,1})^{-1}-z}dH_n(t)+o(1)\\
&=&1+za_{n,0}+o(1).
\eqn
From this, one can easily derive that
\be
\frac1{1+y_na_{n,1}}=1-y_n(1+za_{n,0})+o(1)=1-y_n[1+z\rE m_n(z)]+o(1).
\label{eq3500}
\ee
Finally, from (\ref{eq0013}) with $\ell=0$, we obtain
\be
\rE m_n(z)=\int \frac{1}{t[1-y_n(1+z\rE m_n(z))]-z}dH_n(t)+o(1).
\label{eq353}
\ee
The limiting equation of the above estimate is
\be
m(z)=\int\frac{1}{t[1-y(1+zm_y(z))]-z}dH(t).
\label{eq0014}
\ee
It was proved in \cite{Silv95} that  for each $z\in\mathbb C^+$ the above equation has a unique solution $m(z)$ satisfying $\Im(m)>0$. By this fact, we conclude that $\rE m_n(z)$ tends to the unique solution to the equation (\ref{eq0014}). Finally, by the relation $\um(z)=-(1-y)z^{-1}+ym(z)$, we obtain (\ref{eqg2}).

\subsection{Proof of Theorem \ref{thm2}}

\mgai{First, when the number of columns $k$ in $\bbQ$ is infinite we
  can reduce the situation to the finite $k$ case using arguments
  analoguous to those developed in  Section~\ref{511}.
Therefore, without loss of generality we assume in this section that
$k$ is finite.}

\subsubsection{\mgai{Sketch} of the proof of Theorem \ref{thm2}}
Let $M_n(z)=p[m_n(z)-m_n^0(z)]=n[\um_n(z)-\um_n^0(z)]$
where $m_n(z)$ is the Stieltjes transform of the ESD of $\bbB_n$, $m_n^0(z)$ and $\um_n^0(z)$ are the Stieltjes transforms of $F^{y_n, H_n}$ and $\underline{F}^{y_n, H_n}$ satisfying
$$
z=-\frac1{\um_n^0(z)}+y_n\int\frac{tdH_n(t)}{1+t\um_n^0(z)}
$$
with $y_n=p/n$ and $H_n$ being the ESD of ${\bf T}_n$.
Let $x_r$ be a number greater than
$(1+\sqrt{y})^2\lim\sup\limits_{n}\lambda_{\max}^{\bbT_n}$. Let $x_l$ be a number between $0$ and
$(1-\sqrt{y})^2\lim\inf\limits_{n}\lambda_{\min}^{\bbT_n}$ if the latter is greater than 0 and $y<1$. Otherwise, let $x_l$ be a negative number. Let $\eta_l$ and
$\eta_r$ satisfy
$$x_l<\eta_l<(1-\sqrt{y})^2I_{(0<y<1)}\lim\inf\limits_{n}\lambda_{\min}^{\bbT_n}<
(1+\sqrt{y})^2\lim\sup\limits_{n}\lambda_{\max}^{\bbT_n}<\eta_r<x_r.$$
Let $v_0$ be any positive number. Define a contour ${\cal{C}}=\{x_l+iv: |v|\le
v_0\}\cup{\cal{C}}_u\cup{\cal{C}}_b\cup\{x_r+iv: |v|\le
v_0\}$ where \bqn {\cal{C}}_u=\{x+iv_0: x\in[x_l,x_r]\},\quad
{\cal{C}}_b=\{x-iv_0: x\in[x_l,x_r]\} \eqn and ${\cal C}_n=\{z:
~~z\in{\cal C}~\mbox{and}~|\Im(z)|>n^{-1}\epsilon_n\}$ with
$\epsilon_n\geq n^{-\alpha}$, $0<\alpha<1$. Let
$M_n(z)=p[m_n(z)-m_{F^{y_n, H_n}}(z)]$ and
$$
\widehat{M}_n(z)=\left\{
  \begin{array}{ll}
    M_n(z), & z\in{\cal C}_n,\\
    M_n(x_r+in^{-1}\epsilon_n), & x=x_r,~v\in[0,
    n^{-1}\epsilon_n],\\
    M_n(x_r-in^{-1}\epsilon_n), & x=x_r,~v\in[
    -n^{-1}\epsilon_n,0),\\
    M_n(x_l+in^{-1}\epsilon_n), & x_l<0, x=x_l,~|v|\in[0,
    n^{-1}\epsilon_n],\\\
    M_n(x_l-in^{-1}\epsilon_n), & x_l<0, x=x_l,~|v|\in[-
    n^{-1}\epsilon_n,0).
  \end{array}\right.
$$
\begin{lemma}\label{lem1t}
  Under Assumptions (b)-(c)-(d)-(e)-(f), there exists $\{\epsilon_n\}$ for which
  $\{\widehat M_n(\cdot)\}$ forms a tight sequence on $\cC $.
  Moreover,
  $\widehat M_n(\cdot)$ converges weakly to a two-dimensional Gaussian
  process $M(\cdot)$ satisfying for $z\in\cC $
  \begin{eqnarray*}
    \rE M(z)&=&\alpha_x\frac{y\int{\underline{m}^3(z)t^2}{(1+t\underline{m}(z))^{-3}}dH(t)}
    {\left(1-y\int\frac{\underline{m}^2(z)t^2}{(1+t\underline{m}(z))^{2}}dH(t)\right)\left(1-\alpha_xy\int\frac{\underline{m}^2(z)t^2}
        {(1+t\underline{m}(z))^{2}}dH(t)\right)}\nonumber\\
    &&+\beta_x\frac{y\int{\underline{m}^3(z)t^2}
      {(1+t\underline{m}(z))^{-3}}dH(t)}{1-y\int{\underline{m}^2(z)t^2}
      {(1+t\underline{m}(z))^{-2}}dH(t)},\nonumber
  \end{eqnarray*}
  and variance-covariance function
  \begin{eqnarray*}
    &&{\rm Cov}(M(z_1), M(z_2))\nonumber\\
    &=&\frac{(\partial \underline{m}(z_1)/\partial z_1)(\partial \underline{m}(z_2)/\partial z_2)}{(\underline{m}(z_1)-\underline{m}(z_2))^2}+y\beta_x\left[\int\frac{t}{(\underline{m}(z_1)t+1)^2}
      \frac{t}{(\underline{m}(z_2)t+1)^2}dH(t)\right]\frac{\partial \underline{m}(z_1)}{\partial z_1}\frac{\partial \underline{m}(z_2)}{\partial z_2}\nonumber\\
    &&-\frac{1}{(z_1-z_2)^2}{-}\frac{\partial^2}{\partial z_1\partial z_2}\log(1-a(z_1,z_2))
  \end{eqnarray*}
  where \[a(z_1,z_2)=\alpha_x\left(1+\frac{\underline{m}(z_1)\underline{m}(z_2)(z_1-z_2)}{\underline{m}(z_2)-\underline{m}(z_1)}\right).\]
\end{lemma}
The proof of Lemma \ref{lem1t} is given in the following Section \ref{sec:proofCLT}.
We show how Theorem \ref{thm2} follows from the above Lemma.  We use the
identity
\be\int f(x)dG(x)=-\frac1{2\pi {\bf i}}\oint\limits_{\cC} f(z)m(z)dz\label{tag1.11}\ee
valid for c.d.f. $G$ and $f$ analytic on the support of $G$ where $m(z)$ is the Stieltjes transform of $G$. The complex
integral on the right is over any positively oriented contour enclosing
the support of $G$ and on which $f$ is analytic.   Choose $v_0$, $x_r$, and
$x_l$ so that $f_1,\ldots,f_L$ are all analytic on and inside the
contour $\cC $. Therefore for any $f\in\{f_1,\ldots,f_L\}$, with probability one
$\int f(x)dG_n(x)=-(2\pi {\bf i})^{-1}\int f(z)M_n(z)dz$ for all $n$ large, where the complex integral is over $\cC$.  Moreover, with probability one, for all $n$ large
\bqn\left|\int f(z)(M_n(z)-\widehat M_n(z))dz\right|
&\leq& 4K\epsilon_n(| \max(\lambda_{\max}^{\bbT_n}(1+\sqrt
c_n)^2,\lambda_{\max}^{\bbB_n})-x_r|^{-1} \hfill\\ \hfill
&&+|\min(\lambda_{\min}^{\bbT_n}I_{(0,1)}(c_n)(1-\sqrt c_n)^2,
\lambda_{\min}^{\bbB_n})-x_l|^{-1})\nonumber\eqn
which converges to zero
as $n\to\infty$. Here $K$ is a bound on $f$ over $\cC $. Since
$$\widehat M_n(\cdot)\longrightarrow
\left(-\frac1{2\pi {\bf i}}\int f_1(z)\,\widehat M_n(z)dz,\ldots,
  -\frac1{2\pi {\bf i}}\int f_L(z)\,\widehat M_n(z)dz \right)$$
is a continuous mapping of $C(\cC ,{\mathbb R}^2)$ into ${\mathbb R}^r$, it follows that
the above vector forms tight sequences.  Letting
$M(\cdot)$ denote the limit of any weakly converging subsequence of
$\{\widehat M_n(\cdot)\}$, then we have the weak limit equal in
distribution to
$$\left(-\frac1{2\pi {\bf i}}\int f_1(z)\,M(z)dz,\ldots,
  -\frac1{2\pi {\bf i}}\int f_L(z)\,M(z)dz \right).$$
The fact that this vector,  under the assumptions (b)-(c)-(d)-(e)-(f), is multivariate Gaussian
following from the fact that Riemann sums corresponding to these integrals are
multivariate Gaussian, and that weak limits of Gaussian vectors can only be
Gaussian.   The limiting expressions for the mean and covariance
follow immediately.

\subsubsection{Proof of Lemma \ref{lem1t}}
\label{sec:proofCLT}
The proof of Lemma \ref{lem1t} is similar to the proof of Lemma 1.1 of \cite{BS04} but the proof of the mean and variance-covariance is different from Lemma 1.1 of \cite{BS04}. Let
$M_n(z)=p[m_n(z)-m_n^0(z)]$. In fact, we also have $M_n(z)=n[\um_n(z)-\um_{n}^0(z)]$.

\noindent{\it Convergence of finite dimensional distributions.}
Write for
$z\in\cC _n$, $M_n(z)=M_n^1(z)+M_n^2(z)$ where
$$M_n^1(z)=p[m_n(z)-\rE m_n(z)],\qquad M_n^2(z)=p[\rE m_n(z)-m_{n}^0(z)],$$ and $m_{n}^0(z)$ is the
Stieltjes transform of $F^{y_n,H_n}$. In this subsection we will
show that for any positive integer $r$ and any complex numbers
$z_1,\cdots, z_r$, the random vector $(M_n^{\ell}(z_j),\,j=1,\cdots,r)$
converges to an $2r$-dimensional Gaussian vector for $\ell=1, 2$. Because of
Assumption (e), without loss of generality, we may assume $\|\bbQ\|\leq1$ for all $n$.
Constants appearing in inequalities will be denoted by $K$    and
may take on different values from one expression to the next. Let
\begin{equation}
  \label{eq:epsj}
  \ep_j(z)=\bbr_j^*\bbD^{-1}_j(z)\bbr_j-n^{-1}\rtr\bbT_n\bbD^{-1}_j(z),\quad\delta_j(z)=
  \bbr_j^*\bbD^{-2}_j(z)\bbr_j-n^{-1}\rtr\bbT_n\bbD^{-2}_j(z)=\frac{d}{dz}\ep_j(z),
\end{equation}
and
\begin{equation}
  \label{eq:bn}
  \bar{\beta}_j(z)=\frac{1}{1+n^{-1}\rtr  \bbT_n\bbD_j^{-1}(z)},\quad
  b_n(z)=\frac{1}{1+n^{-1}\rE\rtr  \bbT_n\bbD^{-1}(z)}.
\end{equation}
Notice that
\be
\bbD^{-1}(z)=\bbD_j^{-1}(z)-\bbD_j^{-1}(z)\bbr_j\bbr_j^*\bbD_j^{-1}(z)\beta_j(z).
\label{eqddj}
\ee
By (\ref{eqddj}), we obtain
\begin{eqnarray*}
  p[m_n(z)-\rE m_n(z)]
  &=&\rtr [\bbD^{-1}(z)-\rE \bbD^{-1}(z)]\\
  &=&\sum_{j=1}^n\rtr \rE_j\bbD^{-1}(z)-\rtr \rE_{j-1}\bbD^{-1}(z)\\
  &=&\sum_{j=1}^n\rtr \rE_j[\bbD^{-1}(z)-\bbD^{-1}_j(z)]-
  \rtr \rE_{j-1}[\bbD^{-1}(z)-\bbD^{-1}_j(z)]\\
  &=&-\sum_{j=1}^n(\rE_j-\rE_{j-1})
  \beta_j(z)\bbr_j^*\bbD^{-2}_j(z)\bbr_j\\
  &=&-\frac{d}{dz}\sjln (\rE_j-\rE_{j-1})\log \beta_j(z)\\
  &=&\frac{d}{dz}\sjln (\rE_j-\rE_{j-1})\log \big(1+\ep_j(z)\bar\beta_j(z)\big)
\end{eqnarray*}
where $\beta_j(z)=\bar{\beta}_j(z)-\beta_j(z)\bar{\beta}_j(z)\ep_j(z)$. By Lemma \ref{lem2}, we have
\bqa
&&\rE\left|\frac{d}{dz}\sjln (\rE_j-\rE_{j-1})\big[\log \big(1+\ep_j(z)\bar\beta_j(z)\big)-\ep_j(z)\bar\beta_j(z)\big]\right|^2\non
&\le&\sjln \rE\left|\frac1{2\pi {\bf i}}\oint_{|\zeta-z|=v/2}\frac{\big[\log \big(1+\ep_j(\zeta)\bar\beta_j(\zeta)\big)-\ep_j(\zeta)\bar\beta_j(\zeta)\big]}{(z-\zeta)^2}d\zeta\right|^2\non
&\le&\frac1{2\pi v^4}\sjln \oint_{|\zeta-z|=v/2}\rE|\ep_j(\zeta)\bar\beta_j(\zeta)|^4d\zeta\non
&\le&\frac{K}{2\pi n^4v^4}\sjln \oint_{|\zeta-z|=v/2}\bigg\{\rE\left[\rtr\bbT_n\bbD_j^{-1}(\zeta)\bbT_n\bbD_j^{-1}(z)(\bar \zeta)\right]^2\non
&&+\sum_{i=1}^k \rE|x_{ij}|^8 \rE|\bbq_i^T\bbD_j^{-1}(z)\bbq_i|^4d\zeta\bigg\}\non
&\le&Kn^{-1}+K\eta_n^4\to0.\label{eqtail}
\eqa
Therefore, we need only to derive the finite dimensional limiting distribution of
\be
\frac{d}{dz}\sjln (\rE_j-\rE_{j-1})\ep_j(z)\bar\beta_j(z)=\frac{d}{dz}\sjln \rE_j\ep_j(z)\bar\beta_j(z).
\label{eqgclt}
\ee
Similar to the last three lines of the proof of (\ref{eqtail}), one can show that
\begin{eqnarray}
  &&\sjln\rE\left|(\rE_j-\rE_{j-1})\frac{d}{dz}\ep_j(z)\bar\beta_j(z)\right|^2 I\left(\left|(\rE_j-\rE_{j-1})\frac{d}{dz}\ep_j(z)\bar\beta_j(z)\right|\geq\epsilon\right)\non
  &\leq&\frac{1}{\epsilon^2}\sjln \rE\Big|\rE_j\frac{d}{dz}\ep_j(z)\bar\beta_j(z)\Big|^4\to 0.\nonumber
  \label{eqgclt2}
\end{eqnarray}
Thus, the martingale difference sequence$\{(\rE_j-\rE_{j-1})\frac{d}{dz}\ep_j(z)\bar\beta_j(z)\}$ satisfies the Lyapunov condition. Applying Lemma \ref{lem4},
the random vector $(M_n^1(z_1),\cdots,M_n^1(z_r))$ will tend to an $r$-dimensional Gaussian vector
$(M(z_1),\cdots,$ $ M(z_r))$ whose covariance function is given by
\be
\rCov(M(z_1),M(z_2))=\lim_{n\to\infty}\sjln \rE_{j-1}\left(\rE_j\frac{\partial}{\partial z_1}\ep_j(z_1)\bar\beta_j(z_1)\cdot
  \rE_j\frac{\partial}{\partial z_2}\ep_j(z_2)\bar\beta_j(z_2)\right).\nonumber
\label{eqlimcov}
\ee
Consider the sum
$\Gamma_n(z_1,z_2)=\sum_{j=1}^n\rE_{j-1}[\rE_j(\bar\beta_j(z_1)\epsilon_j(z_1))
\rE_j(\bar\beta_j(z_2)\epsilon_j(z_2))].\label{tag2.6}\nonumber$
Using the same approach of \cite{BS04}, we may replace $\bar\beta_j(z)$ by $b_n(z)$. Therefore, by (1.5) of \cite{BS04}, we have
\bqa\Gamma_n(z_1,z_2)&=&\sum_{j=1}^nb_n(z_1)b_n(z_2)\rE_{j-1}[\rE_j(\epsilon_j(z_1))
\rE_j(\epsilon_j(z_2))]\non
&=&\frac1{n^2}\sum_{j=1}^nb_n(z_1)b_n(z_2)\rtr\Bigg[\rE_{j}(\bbT_n\bbD^{-1}_j(z_1))
\rE_j(\bbT_n\bbD^{-1}_j(z_2))\non
&&\quad\quad\quad\quad\quad\quad\quad\quad\quad+\alpha_x\rtr\rE_j\bar\bbQ\bbQ^*\bbD_j^{-1}(z_1)\bbQ\bbQ^T\rE_j(\bbD_j^T)^{-1}(z_2)\non
&&\quad\quad\quad\quad\quad\quad\quad\quad\quad+\beta_x\sum_{i=1}^k
\bbq_i^*\rE_j\bbD_j^{-1}(z_1)\bbq_i\bbq_i^{*}\rE_j\bbD_j^{-1}(z_2)\bbq_i\Bigg]\non
&=&\frac1{n^2}\sum_{j=1}^nb_n(z_1)b_n(z_2)\rtr\Bigg[\rE_{j}(\bbT_n\bbD^{-1}_j(z_1))
\rE_j(\bbT_n\bbD^{-1}_j(z_2))\non
&&\quad\quad\quad\quad\quad\quad\quad\quad\quad+\alpha_x\rtr\bbT_n\rE_j\bbD_j^{-1}(z_1)\bbT_n\rE_j(\bbD_j^T)^{-1}(z_2)\non
&&\quad\quad\quad\quad\quad\quad\quad\quad\quad+\beta_x\sum_{i=1}^k
\bbq_i^*\rE_j\bbD_j^{-1}(z_1)\bbq_i\bbq_i^{*}\rE_j\bbD_j^{-1}(z_2)\bbq_i\Bigg]
\label{bai2.6}\eqa
where $\alpha_x=|Ex_{11}^2|^2$ and $\beta_x=E|x_{11}^4|-\alpha_x-2$. Here, the third equality holds if either $\alpha_x=0$ or $\bbQ$ is real which implies that $\bbT_n=\bbQ\bbQ^T=\bar\bbQ\bbQ^*$.

Now we use the new method to derive the limit of the first term which is different from but easier than that used in \cite{BS04}.
Let $v_0$ be a lower bound on $\Im(z_i)$. Define $\breve\bbr_j$ as an i.i.d. copy of $\bbr_j$, $j=1,\cdots,n$ and define
$\breve\bbD_j(z)$ similar as $\bbD_j(z)$ by using $\bbr_1,\cdots,\bbr_{j-1},\breve\bbr_{j+1}\cdots,\breve\bbr_n$.
Then we have $$\rtr\Big[\rE_{j}\bbT_n\bbD^{-1}_j(z_1)
\rE_j(\bbT_n\bbD^{-1}_j(z_2))\Big]=\rtr\rE_{j}\Big[\bbT_n\bbD^{-1}_j(z_1)
\bbT_n\breve\bbD^{-1}_j(z_2)\Big].$$
Similar to (\ref{eq3500}), one can prove that
\begin{equation*}
  n^{-1}\rE_j[z_1\rtr\bbT_n\bbD_j^{-1}(z_1)-z_2\rtr\bbT_n\breve\bbD_j^{-1}(z_2)]\to z_1[b^{-1}(z_1)-1]-z_2[b^{-1}(z_2)-1], a.s.
\end{equation*}
where $b(z)=\lim\limits_{n\rightarrow\infty} b_n(z)=-z\um(z)$.
On the other hand,
\begin{eqnarray*}
  &&\rE_jn^{-1}[z_1\rtr\bbT_n\bbD_j^{-1}(z_1)-z_2\rtr\bbT_n\breve\bbD_j^{-1}(z_2)]\\
  &=&n^{-1}\rE_j\rtr\bbT_n\bbD_j^{-1}(z_1)\Big[(z_1-z_2)\sum_{i=1}^{j-1}\bbr_i\bbr_i^*+\sum_{i=j+1}^n(
  z_1\bbr_i\bbr_i^*-z_2\breve\bbr_i\breve\bbr_i^*)\Big]\breve\bbD_j^{-1}(z_2)\\
  &=&n^{-1}\sum_{i=1}^{j-1}(z_1-z_2)\rE_j\bbr_i^*\breve\bbD_{ji}^{-1}(z_2)\bbT_n\bbD_{ji}^{-1}(z_1)\bbr_i\beta_{ji}(z_1)\breve\beta_{ji}(z_2)\\
  &&+n^{-1}\sum_{i=j+1}^n\rE_j\Big[
  z_1\breve\beta_{ji}(z_2)\bbr_i^*\breve\bbD_{ji}^{-1}(z_2)\bbT_n\bbD_j^{-1}(z_1)\bbr_i-
  z_2\beta_{ji}(z_1)\breve\bbr_i^*\breve\bbD_{j}^{-1}(z_2)\bbT_n\bbD_{ji}^{-1}(z_1)\breve\bbr_i\Big]\\
  &=&n^{-2}\sum_{i=1}^{j-1}(z_1-z_2)\rE_j\rtr\bbT_n\breve\bbD_{ji}^{-1}(z_2)\bbT_n\bbD_{ji}^{-1}(z_1)b(z_1)b(z_2)\\
  &&+n^{-2}\sum_{i=j+1}^n\rE_j\Big[
  z_1b(z_2)\rtr\bbT_n\breve\bbD_{ji}^{-1}(z_2)\bbT_n\bbD_j^{-1}(z_1)-
  z_2b(z_1)\rtr\bbT_n\breve\bbD_{j}^{-1}(z_2)\bbT_n\bbD_{ji}^{-1}(z_1)\Big]+o_{a.s.}(1)\\
  &&(\mbox{by replacing}~~\bbD_{ji}^{-1}=\bbD_j^{-1}+\bbD_{ji}^{-1}\bbr_i\bbr_i^*\bbD_{ji}^{-1}\beta_{ji})\\
  &=&\Big[\frac{j-1}n(z_1-z_2)b(z_1)b(z_2)+\frac{n-j}{n}(z_1b(z_2)-z_2b(z_1))\Big]
  \rE_jn^{-1}\rtr\bbT_n\breve\bbD_{j}^{-1}(z_2)\bbT_n\bbD_{j}^{-1}(z_1)+o_{a.s.}(1).
\end{eqnarray*}
Comparing the two estimates, we obtain
$$
\rE_j\frac1n\rtr\bbT_n\breve\bbD_{j}^{-1}(z_2)\bbT_n\bbD_{j}^{-1}(z_1)=\frac{z_1(b^{-1}(z_1)-1)-z_2(b^{-1}(z_2)-1)+o_{a.s.}(1)}
{\frac{j-1}n(z_1-z_2)b(z_1)b(z_2)+\frac{n-j}{n}(z_1b(z_2)-z_2b(z_1))}
$$
Consequently, we obtain
\begin{eqnarray*}
  &&\frac1{n^2}\sum_{j=1}^nb_n(z_1)b_n(z_2)\rtr\rE_{j}\Big[\bbT_n\bbD^{-1}_j(z_1)
  \rE_j(\bbT_n\bbD^{-1}_j(z_2))\Big]\non
  &\to& a(z_1,z_2)\int_0^1\frac1{1-ta(z_1,z_2)}dt={-\log(1-a(z_1,z_2))}=\int_0^{a(z_1,z_2)} \frac1{1-z}dz,
\end{eqnarray*}
where
\begin{eqnarray*}
  a(z_1,z_2)&=&\frac{b(z_1)b(z_2)\bigg[z_1(b^{-1}(z_1)-1)-z_2(b^{-1}(z_2)-1)\bigg]}{z_1b(z_2)
    -z_2b(z_1)}\\
  &=&1+\frac{b(z_1)b(z_2)(z_2-z_1)}{z_1b(z_2)
    -z_2b(z_1)}
  =1+\frac{\um(z_1)\um(z_2)(z_1-z_2)}{\um(z_2)-\um(z_1)}.
\end{eqnarray*}
Thus, we have
\begin{eqnarray*}
  &&\frac{\partial^2}{\partial z_2\partial
    z_1}\int_0^{a(z_1,z_2)} \frac1{1-z}dz\\
  &=&
  \frac{\partial}{\partial z_2}\left(\frac{\frac{\partial}{\partial
        z_1}a(z_1,z_2)}
    {1-a(z_1,z_2)}\right)\\
  &=&\!\frac{\partial}{\partial z_2}\!
  \left[\frac{(\underline{m}(z_2)\!-\!
      \underline{m}(z_1))\bigg(\frac{\partial\underline{m}(z_1)}{\partial z_1}\underline{m}(z_2)(z_1\!-\!z_2)\!+\!
      \underline{m}(z_1)
      \underline{m}(z_2)\bigg)\!+\!\underline{m}(z_1)\underline{m}(z_2)(z_1\!-\!z_2)
      \frac{\underline{m}(z_1)}{\partial z_1}}
    {(\underline{m}(z_2)\!-\!\underline{m}(z_1))^2}\right]\hfill\\
  &&\hfill\shoveright
  {\times\frac{\underline{m}(z_2)-\underline{m}(z_1)}
    {\underline{m}(z_1)\underline{m}(z_2)(z_2-z_1)}}\\
  &=&-\frac{\partial}{\partial z_2}\left(\frac{(\partial\underline{m}(z_1)/\partial z_1)}{\underline{m}(z_1)}+
    \frac1{z_1-z_2}+\frac{(\partial\underline{m}(z_1)/\partial z_1)}{\underline{m}(z_2)-\underline{m}(z_1)}\right)\\
  &=&\frac{(\partial\underline{m}(z_1)/\partial z_1)(\partial\underline{m}(z_2)/\partial z_2)}
  {(\underline{m}(z_2)-\underline{m}(z_1))^2}-\frac1{(z_1-z_2)^2}.
\end{eqnarray*}
Next, we compute the limit of the second term of (\ref{bai2.6}). In this step, we need the assumption that the matrix $\bbQ$ is real.
Similarly, we consider
\begin{eqnarray*}
  &&\frac1n\rE_j[z_1\rtr\bbT_n\bbD_j^{-1}(z_1)-z_2\rtr\bbT_n(\breve\bbD_j^T)^{-1}(z_2)]\to z_1(b^{-1}(z_1)-1)-z_2(b^{-1}(z_2)-1),a.s.
\end{eqnarray*}
On the other hand,
\begin{eqnarray*}
  &&n^{-1}\rE_j[z_1\rtr\bbT\bbD_j^{-1}(z_1)-z_2\rtr\bbT(\breve\bbD_j^T)^{-1}(z_2)]\\
  &=&n^{-1}\rE_j\rtr\bbT\bbD_j^{-1}(z_1)\Big[\sum_{i=1}^{j-1}(z_1\bbr_i\bbr_i^*-z_2\bar\bbr_i\bbr_i^T)+
  \sum_{i=j+1}^n(z_1\bbr_i\bbr_i^*-z_2\bar{\breve\bbr}_i\breve\bbr_i^T)\Big](\breve\bbD_j^T)^{-1}(z_2)\\
  &=&n^{-1}\sum_{i=1}^{j-1}\rE_j\bigg[z_1\breve\beta_{ji}(z_2)\bbr_i^*(\breve\bbD_{ji}^T)^{-1}(z_2)\bbT[\bbD_{ji}^{-1}(z_1)
  -\bbD_{ji}^{-1}(z_1)\bbr_i\bbr_i^*\bbD_{ji}^{-1}(z_1)\beta_{ji}(z_1)]\bbr_i\\
  &&-z_2\bbr_i^T\left((\breve\bbD_{ji}^T)^{-1}(z_2)-\breve\beta_{ji}(z_2)(\breve\bbD^T_{ji})^{-1}(z_2)\bar\bbr_i\bbr_i^T
    (\breve\bbD^T_{ji})^{-1}(z_2)\right)\bbT\bbD_{ji}^{-1}(z_1)\beta_{ji}(z_1)]\bar{\bbr}_i\bigg]\\
  &&+n^{-1}\sum_{i=j+1}^n\rE_j\Big[
  z_1\breve\beta_{ji}(z_2)\bbr_i^*(\breve\bbD_{ji}^T)^{-1}(z_2)\bbT\bbD_j^{-1}(z_1)\bbr_i-
  z_2{\breve\bbr}_i^T\breve\bbD_{j}^{-1}(z_2)\bbT\bbD_{ji}^{-1}(z_1)\bar{\breve\bbr}_i\beta_{ji}(z_1)\Big]\\
  &=&\Big\{(j-1)n^{-1}\alpha_x[-(z_1b(z_2)-z_2b(z_1))+b(z_1)b(z_2)(z_1-z_2)]+(z_1b(z_2)-z_2b(z_1))\Big\}\\
  &&\rE_jn^{-1}\rtr\bbT\bbD_{j}^{-1}(z_1)\bbT(\breve\bbD_{j}^T)^{-1}(z_2)+o_{a.s.}(1).
\end{eqnarray*}
Comparing the two estimates, we obtain
\begin{eqnarray*}
  &&\rE_jn^{-1}\rtr\bbT\breve\bbD_{j}^{-1}(z_2)\bbT\bbD_{j}^{-1}(z_1)\\
  &=&\frac{z_1(b^{-1}(z_1)-1)-(z_2b^{-1}(z_2)-1)+o_{a.s.}(1)}
  {(j-1)n^{-1}\alpha_x[-(z_1b(z_2)-z_2b(z_1))+b(z_1)b(z_2)(z_1-z_2)]+[z_1b(z_2)-z_2b(z_1)]}.
\end{eqnarray*}
Consequently, we obtain
\begin{eqnarray*}
  &&n^{-2}\sum_{j=1}^n\alpha_xb_n(z_1)b_n(z_2)\rtr\rE_{j}\bbT\bbD^{-1}_j(z_1)
  \rE_j(\bbT(\bbD^T_j)^{-1}(z_2))\non
  &\to& \tilde a(z_1,z_2)\int_0^1\frac1{1-t\tilde a(z_1,z_2)}dt={-}\log(1-\tilde a(z_1,z_2))=\int\limits_0^{\tilde a(z_1,z_2)}\frac{1}{1-z}dz,
\end{eqnarray*}
where
\begin{eqnarray*}
  \tilde a(z_1,z_2)&=&\frac{\alpha_xb(z_1)b(z_2)(z_1(b^{-1}(z_1)-1)-z_2(b^{-1}(z_2)-1))}{z_1b(z_2)
    -z_2b(z_1)}\\
  &=&\alpha_x\left(1+\frac{b(z_1)b(z_2)(z_2-z_1)}{z_1b(z_2)
      -z_2b(z_1)}\right)
  =\alpha_x\left(1+\frac{\um(z_1)\um(z_2)(z_1-z_2)}{\um(z_2)-\um(z_1)}\right).
\end{eqnarray*}
Last, we will compute the limit of the third term of (\ref{bai2.6}). By (9.9.12) of \cite{BS10}, we have
\begin{eqnarray*}
  &&n^{-2}\sum_{j=1}^{n}\beta_x\sum\limits_{i=1}^{k}{\bf e}_i^T\bbQ^*\bbD_j^{-1}(z_1)\bbQ{\bf e}_i{\bf e}_i^T\bbQ^*\bbD_j^{-1}(z_2)\bbQ{\bf e}_i\\
  &=&\frac{1}{n^2z_1z_2}\sum_{j=1}^{n}\beta_x\sum\limits_{i=1}^{k}{\bf e}_i^T\bbQ^*(\underline{m}(z_1)\bbT_n+\bbI_p)^{-1}\bbQ{\bf e}_i{\bf e}_i^T\bbQ^*(\underline{m}(z_2)\bbT_n+\bbI_p)^{-1}\bbQ{\bf e}_i+o_p(1)\\
\end{eqnarray*}
If $\beta_x\ne 0$, then
\mgai{by assumption, the matrix $\bbQ^{*}\bbQ$ is diagonal.}
Using the identity
$\bbQ^*[\underline{m}(z)\bbT_n+\bbI_p]^{-1}\bbQ=\underline{m}^{-1}(z)\{\bbI_k-[\underline{m}(z)\bbQ^*\bbQ+\bbI_k]^{-1}\}$,
\mgai{we see that  the matrix is diagonal. We have}
\begin{eqnarray*}
  &&\frac{1}{n^2z_1z_2}\sum_{j=1}^{n}\beta_x\sum\limits_{i=1}^{k}{\bf e}_i^T\bbQ^*[\underline{m}(z_1)\bbT_n+\bbI_p]^{-1}\bbQ{\bf e}_i{\bf e}_i^T\bbQ^*[\underline{m}(z_2)\bbT_n+\bbI_p]^{-1}\bbQ{\bf e}_i\\
  &=&\frac{1}{n^2z_1z_2}\sum_{j=1}^{n}\beta_x{\rm tr}\{\bbQ^*[\underline{m}(z_1)\bbT_n+\bbI_p]^{-1}\bbQ\bbQ^*[\underline{m}(z_2)\bbT_n+\bbI_p]^{-1}\bbQ\}\\
  &=&\frac{y\beta_x}{z_1z_2}\int\frac{t^2}{[1+t\underline{m}(z_1)][1+t\underline{m}(z_2)]}dH(t)+o(1).
\end{eqnarray*}
Then the third term of $\Cov(M(z_1), M(z_2))$ is
\begin{eqnarray*}
  &&\frac{\partial^2}{\partial z_1\partial z_2}\left\{\frac{y\beta_xb(z_1)b(z_2)}{z_1z_2}\int\frac{t^2}{[1+t\underline{m}(z_1)][1+t\underline{m}(z_2)]}dH(t)\right\}\\
  &=&\frac{\partial^2}{\partial z_1\partial z_2}\left\{y\beta_x\int\frac{t^2\underline{m}(z_1)\underline{m}(z_2)}{[1+t\underline{m}(z_1)][1+t\underline{m}(z_2)]}dH(t)\right\}\\
  &=&y\beta_x\frac{\partial\underline{m}(z_1)}{\partial z_1}
  \frac{\partial\underline{m}(z_2)}{\partial z_2}\int\frac{t^2}{[1+t\underline{m}(z_1)]^2[1+t\underline{m}(z_2)]^2}dH(t).
\end{eqnarray*}

{\it Tightness of $M_n^1(z)$.}
As done in \cite{BS04}, the proof of tightness of $M_n^1$ relies on the proof of
\bqa
&& \sup_{n,z_1,z_2\in\cC}\rE\left|\frac{\hat M_n^1(z_1)-\hat M_n(z_2)-\rE(\hat M_n(z_1)-\hat M_n(z_2))}{z_1-z_2}\right|^2\non
&=&\sup_{n,z_1,z_2\in\cC\atop |\Im(z_i)|>n^{-1}\epsilon_n}\rE\left|\rtr \bbD^{-1}(z_1)\bbD^{-1}(z_2)-\rE\rtr \bbD^{-1}(z_1)\bbD^{-1}(z_2)\right|^2\non
&:=&\sup_{n,z_1,z_2\in\cC\atop |\Im(z_i)|>n^{-1}\epsilon_n} J_n(z_1,z_2)<\infty
\label{tightg}
\eqa
where by the formula (\ref{eqddj}), we have
\begin{eqnarray*}
  &&J_n(z_1,z_2)\\
  &=&\rE\left|\sum_{j=1}^n(\rE_j-\rE_{j-1}){\rm tr}\bbD^{-1}(z_1)\bbD^{-1}(z_2)\right|^2\\
  &=&\sum_{j=1}^n\rE\left|(\rE_j-\rE_{j-1})\bbr_j^*\bbD_j^{-1}(z_1)\bbD_j^{-1}(z_2)\bbD_j^{-1}(z_1)\bbr_j\beta_j(z_1)\right|^2\\
  &&+\sum_{j=1}^n\rE\left|(\rE_j-\rE_{j-1})\bbr_j^*\bbD_j^{-1}(z_2)\bbD_j^{-1}(z_1)\bbD_j^{-1}(z_2)\bbr_j\beta_j(z_2)\right|^2\\
  &&+\sum_{j=1}^n\rE\left|(\rE_j-\rE_{j-1})\bbr_j^*\bbD_j^{-1}(z_1)\bbD_j^{-1}(z_2)\bbr_j\bbr_j^*\bbD_j^{-1}(z_2)\bbD_j^{-1}(z_1)\bbr_j\beta_j(z_1)\beta_j(z_2)\right|^2\\
  &=&\sum_{j=1}^n\rE\left|(\rE_j-\rE_{j-1})\bbr_j^*\bbD_j^{-1}(z_1)\bbD_j^{-1}(z_2)\bbD_j^{-1}(z_1)\bbr_jb_j(z_1)\right|^2\\
  &&+\sum_{j=1}^n\rE\left|(\rE_j-\rE_{j-1})\bbr_j^*\bbD_j^{-1}(z_2)\bbD_j^{-1}(z_1)\bbD_j^{-1}(z_2)\bbr_jb_j(z_2)\right|^2\\
  &&+\sum_{j=1}^n\rE\left|(\rE_j-\rE_{j-1})\bbr_j^*\bbD_j^{-1}(z_1)\bbD_j^{-1}(z_2)\bbr_j\bbr_j^*\bbD_j^{-1}(z_2)\bbD_j^{-1}(z_1)\bbr_jb_j(z_1)b_j(z_2)\right|^2+O(1)\\
  &=&O(1).
\end{eqnarray*}

{\it Convergence of $M_n^2(z)$.}
In order to simplify the exposition,
we let $\cC _1=\cC _u$ or $\cC _u\cup\cC _l$ if $x_l<0$, and $\cC _2=\cC _r$
or $\cC _r\cup\cC _l$ if $x_l>0$.  We begin with proving
\be\sup_{z\in\cC _n}|\rE \underline{m}_n(z)-\underline{m}(z)|\to0\quad\text{ as $n\to\infty$.}
\label{tag4.1}\ee
Since $F^{\underline \bbB_n}\darrow F^{y,H}$ almost surely, we get from
d.c.t. $\rE F^{\underline \bbB_n}\darrow F^{y,H}$. It is easy to
verify that $\rE F^{\underline \bbB_n}$ is a proper c.d.f.  Since, as $z$
ranges in $\cC _1$, the functions $(\lambda-z)^{-1}$ in
$\lambda\in[0,\infty)$ form a bounded, equicontinuous family, it follows
[see, e.g. Billingsley (1968), Problem 8, p. 17] that
$$\sup_{z\in\cC _1}|\rE\underline{m}_n(z)-\underline{m}(z)|\to0.$$
For $z\in\cC _2$, we write ($\eta_l,\eta_r$ defined as in previous section)
$$\rE m_n(z)-m(z)=\int\frac1{\lambda-z}I_{[\eta_l,\eta_r]}(\lambda)
d(\rE F^{\underline \bbB_n}(\lambda)-F^{c,H}(\lambda))+
\rE\int\frac1{\lambda-z}I_{[\eta_l,\eta_r]^c}(\lambda)dF^{\underline \bbB_n}
(\lambda).$$
As above, the first term converges uniformly to zero.  With $\ell\ge2$, we get
\begin{eqnarray*}
  &&\sup_{z\in\cC _2}\left|\rE\int\frac1{\lambda-z}I_{[\eta_l,\eta_r]^c}(\lambda)
    dF^{\underline \bbB_n}(\lambda)\right|\hfill\\ &&\hfill
  \leq (\epsilon_n/n)^{-1}P(\|\bbB_n\|\ge\eta_r\text{ or }
  \lambda_{\min}^{\bbB_n}\leq\eta_l)
  \leq Kn\epsilon_n^{-1}n^{-\ell} \to0.
\end{eqnarray*}
From the fact that $F^{y_n,H_n}\darrow  F^{y,H}$ [see \cite{BS98}, below (3.10)] along with the fact that $\cC $ lies outside the support of $F^{c,H}$,
it is straightforward to verify that
\be\sup_{z\in\cC }|{\underline{m}_n^0(z)}-\underline{m}(z)|\to0\quad\text{ as $n\to\infty$}
\label{tag4.2}\ee
where $m_n^0(z)$ is the Stieltjes transform of $F^{y_n,H_n}$.
We now show that
\be\sup_{n, z\in\cC _n}\|({\rE\mgai{\underline{m}_n}(z)}\bbT_n+\bbI_p)^{-1}\|<\infty.\label{tag4.3}
\ee
From
Lemma 2.11 of \cite{BS98}, $\|({\rE\mgai{\underline{m}_n}(z)}\bbT_n+\bbI_p)^{-1}\|$ is bounded by $\max(2,4v_0^{-1}$) on $\cC _u$.
On $\cC _l$, the boundedness of $\|({\rE\mgai{\underline{m}_n}(z)}\bbT_n+\bbI_p)^{-1}\|$ follows from the fact that
$$
|\rE\umn(z)\lambda_j+1|\ge \Re(\rE\umn(z)\lambda_j+1)\ge P(\lambda_{\min}(\bbB_n)>\eta_l)-\epsilon_n^{-1}P(\lambda_{\min}(\bbB_n)<\eta_l)
\to 1,
$$
where $\lambda_j$ is an arbitrary eigenvalue of $\bT_n$ {and $\|\bbQ\|\leq 1$}. Now let us consider the bound on $\cC _r$. By (1.1) of
\cite{BS98}, there exists a support point $t_0$ of $H$ such that
$1+t_0\underline{m}(z)\ne 0$ on $\cC _r$. Since $\underline{m}(z)$ is analytic on $\cC _r$, there exist positive constants
$\delta_1$ and $\mu_0$ such that
$$
\inf_{z\in \cC _r}|1+t_0\underline{m}(z)|>\delta_1,\ \ \text{ and }\ \ \sup_{z\in\cC _r}|\underline{m}(z)|<\mu_0.
$$
By (\ref{tag4.1}) and $H_n\to H$, for all large $n$, there exists an integer $j\le n$ such that
$|\lambda_j-t_0|<\delta_1/4\mu_0$ and $\sup_{z\in\cC _r}|E\underline{m}_n(z)-\underline{m}(z)|<\delta_1/4$. Then, we have
$$
\inf_{z\in \cC _r}|1+\lambda_j\rE\underline{m}_n(z)|>\delta_1/2,
$$
which completes the proof of (\ref{tag4.3}).

Next we show the existence of $\xi\in(0,1)$ such that for all $n$ large
\be \sup_{z\in\cC _n}
\left|y_n\rE\umn(z)^2\int\frac{t^2}{(1+t\rE\umn(z))^2}dH_n(t)
\right|<\xi.\label{tag4.4}\ee
From the identity (1.1) of Bai and Silverstein (1998)
$$\um(z)=\frac1{-z+y\int\frac{t}{1+t\um(z)}dH(t)}$$
valid for $z=x+iv$ outside the support of $F^{y,H}$, we find
$$\Im\,\um(z)=\frac{v+\Im\,\um(z)y\int\frac{t^2}{|1+t\um(z)|^2}dH(t)}
{\left|-z+y\int\frac{t}{1+t\um(z)}dH(t)\right|^2}.$$
Therefore
\begin{eqnarray}
  \left|y\um(z)^2\int\frac{t^2}{(1+t\um(z))^2}dH(t)\right|
  &\leq&\frac{y\int\frac{t^2}{|1+t\um(z)|^2}dH(t)}
  {\left|-z+y\int\frac{t}{1+t\um(z)}dH(t)\right|^2}\non
  &=&\frac{\Im\,\um(z)y\int\frac{t^2}{|1+t\um(z)|^2}dH(t)}
  {v+\Im\,\um(z)y\int\frac{t^2}{|1+t\um(z)|^2}dH(t)}\non
  &=&\frac{y\int\frac{1}{|x-z|^2}dF^{y,H}(x)\int\frac{t^2}{|1+t\um(z)|^2}dH(t)}
  {1+y\int\frac{1}{|x-z|^2}dF^{y,H}(x)\int\frac{t^2}{|1+t\um(z)|^2}dH(t)}<1,\label{tag4.5}
\end{eqnarray}
for all $z\in \cC $. By continuity, we have $\xi_1<1$ such that
\be \sup_{z\in \cC }\left|y\um(z)^2\int\frac{t^2}{(1+t\um(z))^2}dH(t)\right|<\xi_1.\label{tag4.6}\ee
Therefore, using (\ref{tag4.1}), (\ref{tag4.4}) follows.
We proceed with some improved bounds on quantities appearing earlier.
\mgai{Recall  the functions $\beta_j(z)$ and $\gamma_j(z)$ defined in \eqref{eq:betaj-gammaj}.}
For $p\geq4$, we have
\be\rE|\gamma_1(z)|^p\leq Kn^{-2}.\label{tag4.7}\ee
Let $\bbM$ be nonrandom $p\times p$.  Then
\begin{eqnarray*}
  &&\rE|\rtr\bbD^{-1}(z)\bbM-\rE\rtr\bbD^{-1}(z)\bbM|^2\\
  &=&\rE|\sum_{j=1}^n\rE_j\rtr\bbD^{-1}(z)\bbM-\rE_{j-1}\rtr\bbD^{-1}(z)\bbM|^2\\
  &=&\rE|\sum_{j=1}^n(\rE_j-\rE_{j-1})\rtr(\bbD^{-1}(z)-\bbD_j^{-1}(z))\bbM|^2\\
  &=&\sum_{j=1}^n\rE|(\rE_j-\rE_{j-1})\beta_j(z)\bbr_j^*\bbD_j^{-1}(z)\bbM\bbD_j^{-1}(z)\bbr_j|^2\\
  &=&|b_n(z)|^2\sum_{j=1}^n\rE|(\rE_j-\rE_{j-1})(1-\beta_j(z)\gamma_j(z))\bbr_j^*\bbD_j^{-1}(z)\bbM\bbD_j^{-1}(z)\bbr_j|^2\\
  &=&|b_n(z)|^2\sum_{j=1}^n\rE|\rE_j(\bbr_j^*\bbD_j^{-1}(z)\bbM\bbD_j^{-1}(z)\bbr_j-n^{-1}\rtr\bbD_j^{-1}(z)\bbM\bbD_j^{-1}(z)\bbT_n\\
  &&-(\rE_j-\rE_{j-1})\beta_j(z)\gamma_j(z)\bbr_j^*\bbD_j^{-1}(z)\bbM\bbD_j^{-1}(z)\bbr_j|^2\\
  &\leq&2|b_n(z)|^2n\rE|\bbr_j^*\bbD_j^{-1}(z)\bbM\bbD_j^{-1}(z)\bbr_j-n^{-1}\rtr\bbD_j^{-1}(z)\bbM\bbD_j^{-1}(z)\bbT_n|^2\\
  &&+4|b_n(z)|^2n\rE|\beta_j(z)\gamma_j(z)\bbr_j^*\bbD_j^{-1}(z)\bbM\bbD_j^{-1}(z)\bbr_j|^2\\
  &\leq&2|b_n(z)|^2n\rE|\bbr_j^*\bbD_j^{-1}(z)\bbM\bbD_j^{-1}(z)\bbr_j-n^{-1}\rtr\bbD_j^{-1}(z)\bbM\bbD_j^{-1}(z)\bbT_n|^2\\
  && +4|b_n(z)|^2n(\rE|\gamma_j(z)|^4)^{1/2}(\rE|\beta_j(z)|^8)^{1/4}(\rE|\bbr_j^*\bbD_j^{-1}(z)\bbM\bbD_j^{-1}(z)\bbr_j|^8)^{1/4}.
\end{eqnarray*}
Using (9.9.3) of \cite{BS10}, (\ref{tag4.7}) and the boundedness of $b_n(z)$ we get
\be\rE|\rtr\bbD^{-1}(z)\bbM-\rE\rtr\bbD^{-1}(z)\bbM|^2\leq K\|\bbM\|^2.\label{tag4.8}\ee
The same argument holds for $\bbD_1^{-1}(z)$ so we also have
\be\rE|\rtr\bbD_1^{-1}(z)\bbM-\rE\rtr\bbD_1^{-1}(z)\bbM|^2\leq K\|\bbM\|^2.\label{tag4.9}\ee
Our next task is to investigate the limiting behavior of
\begin{eqnarray*}
  && n\left(y_n\int\frac{dH_n(t)}{1+t\rE\umn(z)}+zy_n\rE m_n(z)\right)\\
  &=&{n\rE\beta_1(z)\biggl[\bbr_1^*\bbD_1^{-1}(z)(\rE\umn(z)\bbT_n+\bbI_p)^{-1}\bbr_1
    -n^{-1}\rE\rtr(\rE\umn(z)\bbT_n+\bbI_p)^{-1}\bbT_n\bbD^{-1}(z)\biggr]},
\end{eqnarray*}
for $z\in\cC _n$ [see (5.2) in \cite{BS98}]. Throughout
the following, all bounds, including $O(\cdot)$ and $o(\cdot)$ expressions,
and convergence statements hold uniformly for $z\in\cC _n$.
We have
\begin{eqnarray}
  &&\rE\rtr(\rE\umn(z)\bbT_n+\bbI_p)^{-1}\bbT_n\bbD_1^{-1}(z)-\rE\rtr(\rE\umn(z)\bbT_n+\bbI_p)^{-1}\bbT_n\bbD^{-1}(z)\label{tag4.10}\\
  &=&\rE\beta_1(z)\rtr(\rE\umn(z)\bbT_n+\bbI_p)^{-1}\bbT_n\bbD_1^{-1}(z)\bbr_1\bbr_1^*\bbD_1^{-1}(z)\nonumber\\
  &=&b_n(z)\rE(1-\beta_1(z)\gamma_1(z))\bbr_1^*\bbD_1^{-1}(z)(\rE\umn(z)\bbT_n+\bbI_p)^{-1}\bbT_n\bbD_1^{-1}(z)\bbr_1.\nonumber
\end{eqnarray}
From (\ref{tag4.3}) we get
\begin{eqnarray*}
  &&|\rE\beta_1(z)\gamma_1(z)\bbr_1^*\bbD_1^{-1}(z)(\rE\umn(z)\bbT_n+\bbI_p)^{-1}
  \bbT_n\bbD_1^{-1}(z)\bbr_1|\\ &\leq&(\rE|\gamma_1(z)|^2)^{1/2})(\rE|\beta_1(z)|^4)^{1/4})
  (\rE|\bbr_1^*\bbD_1^{-1}(z)(\rE\umn(z)\bbT_n+\bbI_p)^{-1}\bbT_n\bbD_1^{-1}(z)\bbr_1|^4)^{1/4}\\
  &\leq& Kn^{-1/2}.
\end{eqnarray*}
Therefore
$$|(\ref{tag4.10})-n^{-1}b_n(z)\rE\rtr\bbD_1^{-1}(z)(\rE\umn(z)\bbT_n+\bbI_p)^{-1}\bbT_n\bbD_1^{-1}(z)\bbT_n|\leq  Kn^{-1/2}.$$
Since $\beta_1(z)=b_n(z)-b_n^2(z)\gamma_1(z)+\beta_1(z)b_n^2(z)\gamma_1^2(z)$, we have
\begin{eqnarray*}
  &&n\rE\beta_1(z)\bbr_1^*\bbD_1^{-1}(z)(\rE\umn(z)\bbT_n+\bbI_p)^{-1}\bbr_1-
  \rE\beta_1(z)\rE\rtr(\rE\umn(z)\bbT_n+\bbI_p)^{-1}\bbT_n\bbD_1^{-1}(z)\\
  &=&-b_n^2(z)n\rE\gamma_1(z)\bbr_1^*\bbD_1^{-1}(z)(\rE\umn(z)\bbT_n+\bbI_p)^{-1}\bbr_1\\
  &&+b_n^2(z)(n\rE\beta_1(z)\gamma_1^2(z)\bbr_1^*\bbD_1^{-1}(z)(\rE\umn(z)\bbT_n+\bbI_p)^{-1}\bbr_1\\
  &&-(\rE\beta_1(z)\gamma_1^2(z))\rE\rtr(\rE\umn(z)\bbT_n+\bbI_p)^{-1}\bbT_n\bbD_1^{-1}(z))\\
  &=&-b_n^2(z)n\rE\gamma_1(z)\bbr_1^*\bbD_1^{-1}(z)(\rE\umn(z)\bbT_n+\bbI_p)^{-1}\bbr_1\\
  &&+b_n^2(z)(\rE[n\beta_1(z)\gamma_1^2(z)\bbr_1^*\bbD_1^{-1}(z)(\rE\umn(z)\bbT_n+\bbI_p)^{-1}\bbr_1\\
  &&-\beta_1(z)\gamma_1^2(z)\rtr\bbD_1^{-1}(z)(\rE\umn(z)\bbT_n+\bbI_p)^{-1}\bbT_n])\\
  &&+b_n^2(z)\text{Cov}(\beta_1(z)\gamma_1^2(z),\rtr\bbD_1^{-1}(z)(\rE\umn(z)\bbT_n+\bbI_p)^{-1}\bbT_n)
\end{eqnarray*}
($\text{Cov}(X,Y)=\rE XY-\rE X\rE Y$). Using (\ref{tag4.3}), (\ref{tag4.7}), (\ref{tag4.9}) and Lemma \ref{lem2}, we have
\begin{eqnarray*}
  &&|\rE[n\beta_1(z)\gamma_1^2(z)\bbr_1^*\bbD_1^{-1}(z)(\rE\umn(z)\bbT_n+\bbI_p)^{-1}\bbr_1-\beta_1(z)\gamma_1^2(z)\rtr\bbD_1^{-1}(z)(\rE\umn(z)\bbT_n+\bbI_p)^{-1}\bbT_p]|\\
  &\leq& n(\rE|\beta_1(z)|^4)^{1/4}(\rE|\gamma_1(z)|^4)^{1/2}(\rE|\bbr_1^*\bbD_1^{-1}(z)(\rE\umn(z)\bbT_n+\bbI_p)^{-1}\bbr_1\\
  &&-n^{-1}\rtr\bbD_1^{-1}(\rE\umn(z)\bbT_n+\bbI_p)^{-1}\bbT_n|^4)^{1/4}\\
  &\leq& Kn^{-1/2}
\end{eqnarray*}
and
\begin{eqnarray*}
  &&|\text{Cov}(\beta_1(z)\gamma_1^2(z),\rtr\bbD_1^{-1}(z)(\rE\umn(z)\bbT_n+\bbI_p)^{-1}\bbT_p)|\\
  &\leq&(\rE|\beta_1(z)|^6)^{1/6}(\rE|\gamma_1(z)|^6)^{1/3}(\rE|\rtr\bbD_1^{-1}(z)(\rE\umn(z)\bbT_n+\bbI_p)^{-1}\bbT_n\\
  &&-\rE\rtr\bbD_1^{-1}(z)(\rE\umn(z)\bbT_n+\bbI_p)^{-1}\bbT_n|^2)^{1/2}\\
  &\leq& Kn^{-2/3}.
\end{eqnarray*}
Since $\beta_1(z)=b_n(z)-b_n(z)\beta_1(z)\gamma_1(z)$, then we have $\rE\beta_1(z)=b_n(z)+O(n^{-1/2})$.
Write \begin{eqnarray*}
  &&\rE n\gamma_1(z)\bbr_1^*\bbD_1^{-1}(z)(\rE\umn(z)\bbT_n+\bbI_p)^{-1}\bbr_1\\
  &=&n\rE[(\bbr_1^*\bbD_1^{-1}(z)\bbr_1-n^{-1}\rtr\bbD_1^{-1}(z)\bbT_n)(\bbr_1^*\bbD_1^{-1}(z)(\rE\umn(z)\bbT_n+\bbI_p)^{-1}\bbr_1\\
  &&-n^{-1}\rtr\bbD_1^{-1}(z)(\rE\umn(z)\bbT_n+\bbI_p)^{-1}\bbT_n)]\\
  &&+n^{-1}\text{Cov}(\rtr\bbD_1^{-1}(z)\bbT_n,\rtr\bbD_1^{-1}(z)(\rE\umn(z)\bbT_n+\bbI_p)^{-1}\bbT_n).
\end{eqnarray*} From (\ref{tag4.9}) we see that the second
term above is $O(n^{-1})$.
Therefore, we arrive at
\begin{eqnarray*}
  &&n\left(y_n\int\frac{dH_n(t)}{1+t\rE\umn}+zy_n\rE m_n(z)\right)\\
  &=&b_n^2(z)n^{-1}\rE\rtr\bbD_1^{-1}(z)(\rE\umn(z)\bbT_n+\bbI_p)^{-1}\bbT_n\bbD_1^{-1}(z)\bbT_n\hfill\non
  &&-b_n^2(z)n\rE[(\bbr_1^*\bbD_1^{-1}(z)\bbr_1-n^{-1}\rtr\bbD_1^{-1}(z)\bbT_n)(\bbr_1^*\bbD_1^{-1}(z)(\rE\umn(z)\bbT_n+\bbI_p)^{-1}\bbr_1\\
  &&-n^{-1}\rtr\bbD_1^{-1}(z)(\rE\umn(z)\bbT_n+\bbI_p)^{-1}\bbT_n)]+o(1)\non
  &=&-\frac{b_n^2(z)}{n}\beta_x\sum\limits_{i=1}^{m}\bbq_i^{*}\bbD_1^{-1}(z)\bbq_i\cdot\bbq_i^{*}\bbD_1^{-1}(z)(\rE\underline{m}_n(z)\bbT_n+\bbI)^{-1}\bbq_i\label{m1}\\
  &&-\frac{b_n^2(z)}{n}\alpha_x\rtr\bbT_n(\bbD_1^T(z))^{-1}\bbT_n\bbD_1^{-1}(z)(\rE\underline{m}_n(z)\bbT_n+\bbI)^{-1}.
\end{eqnarray*}
Let $A_n(z)=y_n\int\frac{dH_n(t)}{1+t\rE\umn(z)}+zy_n\rE m_n(z)$. Using the identity
$\rE\umn(z)=-\frac{(1-y_n)}z+y_n\rE m_n(z)$, we have
\begin{eqnarray*}
  A_n(z)&=&y_n\int\frac{dH_n(t)}{1+t\rE\umn(z)}-y_n+z\rE\umn(z)+1
  \hfill\\ && \hfill=-\rE\umn(z)\left(-z-\frac1{\rE\umn(z)}+
    y_n\int\frac{tdH_n(t)}{1+t\rE\umn(z)}\right).
\end{eqnarray*}
It follows that
$$\rE\umn(z)=
\frac1{-z+y_n\int\frac{tdH_n(t)}{1+t\rE\umn(z)}+A_n/\rE\umn(z)}.$$
From this, together with the analogous identity [below (4.4)]
we get
\be\rE\umn(z)-\umn^0(z)=-\frac{\umn^0(z)A_n}{1-y_n\rE\umn(z)\umn^0(z)\int\frac{t^2dH_n(t)}
  {(1+t\rE\umn(z))(1+t\umn^0(z))}}.\label{tag4.12}\ee
We see from (\ref{tag4.4}) and the corresponding bound involving $\um_n^0(z)$,
that the denominator of (\ref{tag4.12}) is bounded away from zero.
First we have
\begin{eqnarray*}
  &&\alpha_x\frac{z^2\underline{m}^2(z)}{n^2}\sum\limits_{j=1}^{n}\rE\rtr[\bbD_1^{-1}(z)\{\underline{m}(z)\bbT_n+{\bf I}_p\}^{-1}\bbT_n(\bbD_1^T(z))^{-1}\bbT_n]\\
  &=&\alpha_x\frac{\underline{m}^2(z)}{n}\rE\rtr(\underline{m}(z)\bbT_n+{\bf I}_p)^{-3}\bbT_n^2\\
  &&+\alpha_x\frac{z^2\underline{m}^4(z)}{n^2}\sum\limits_{i\not=j}
  \rE\rtr\bigg\{(\underline{m}(z)\bbT_n+{\bf I}_p)^{-1}\bbT_n\underline{m}(z)\bbT_n+{\bf I}_p)^{-1}(\bbr_i\bbr_i^{*}-\frac{1}{n}\bbT_n)\bbA_{ij}^{-1}(z)\\
  &&\quad\quad\quad\quad\quad\quad\quad\quad\quad\cdot(\underline{m}(z)\bbT_n+{\bf I}_p)^{-1}\bbT_n(\bbA^T_{ij})^{-1}(z)(\bar\bbr_i\bbr_i^T-\frac{1}{n}\bbT_n)\bigg\}+o(1)\\
  &=&\alpha_x\frac{\underline{m}^2(z)}{n}\rtr(\underline{m}(z)\bbT_n+{\bf I}_p)^{-3}\bbT_n^2\\
  &&+\frac{\alpha_x^2z^2\underline{m}^4(z)}{n^3}\sum\limits_{j=1}^{n}\{\rtr(\underline{m}(z)\bbT_n+{\bf I}_p)^{-1}\bbT_n\}^2\cdot\rtr\bbD_1^{-1}(z)(\underline{m}(z)\bbT_n+{\bf I}_p)^{-1}\bbT_n(\bbA^T_{j}(z))^{-1}\bbT_n+o(1).
\end{eqnarray*}
Then we have
\begin{eqnarray*}
  &&\alpha_x\frac{z^2\underline{m}^2(z)}{n^2}\sum\limits_{j=1}^{n}\rE\rtr[\bbD_1^{-1}(z)\{\underline{m}(z)\bbT_n+{\bf I}_p\}^{-1}\bbT_n(\bbD_1^T(z))^{-1}\bbT_n]\\
  &=&\frac{\alpha_x\frac{\underline{m}^2(z)}{n}\rtr\{(\underline{m}(z)\bbT_n+{\bf I}_p)^{-3}\bbT_n^2\}}{1-\frac{\alpha_x\underline{m}^2(z)}{n}\rtr\{(\underline{m}(z)\bbT_n+{\bf I}_p)^{-1}\bbT_n\}^2}+o(1)\\
  &=&\frac{\alpha_xy\int\frac{\underline{m}^2(z)t^2dH(t)}{(1+t\underline{m}(z))^3}}{1-\alpha_xy\int\frac{\underline{m}^2(z)t^2dH(t)}{(1+t\underline{m}(z))^2}}+o(1).
\end{eqnarray*}
Thus we obtain
\begin{eqnarray*}
  &&\frac{\alpha_x\frac{z^2\underline{m}^2(z)}{n^2}\sum\limits_{j=1}^{n}\rE\rtr[\bbD_1^{-1}(z)\{\underline{m}(z)\bbT_n+{\bf I}_p\}^{-1}\bbT_n(\bbD_1^T(z))^{-1}\bbT_n]}{1-y\int\frac{\underline{m}^2(z)t^2dH(t)}{(1+t\underline{m}(z))^2}}\\
  &=&\frac{\alpha_xy\int\frac{\underline{m}^2(z)t^2dH(t)}{(1+t\underline{m}(z))^3}}{\{1-y\int\frac{\underline{m}^2(z)t^2dH(t)}
    {(1+t\underline{m}(z))^2}\}\{1-\alpha_xy\int\frac{\underline{m}^2(z)t^2dH(t)}{(1+t\underline{m}(z))^2}\}}+o(1).
\end{eqnarray*}
{Moreover, if $\beta_x=0$ or $\bbQ^*{\bbQ}$ is diagonal, then we have
  \begin{eqnarray*}
    &&\frac{\beta_xy_nz^2\underline{m}^3(z)\frac{1}{pn}\sum\limits_{j=1}^{n}\sum\limits_{i=1}^{m}
      \rE\bbq_i^{*}\bbD_1^{-1}(z)\{\underline{m}(z)\bbT_n+{\bf I}_p\}^{-1}\bbq_i\cdot\bbq_i^{*}\bbD_1^{-1}(z)\bbq_i}{1-y\int\frac{\underline{m}^2(z)t^2dH(t)}{(1+t\underline{m}(z))^2}}\\
    &=&\beta_x \frac{y\int\frac{\underline{m}^3(z)t^2}
      {(\underline{m}(z)t+1)^3}dH(t)}{1-y\int\frac{\underline{m}^2(z)t^2dH(t)}
      {(1+t\underline{m}(z))^2}}+o(1)
  \end{eqnarray*}}
where
\begin{eqnarray*}
  &&\bbq_i^{*}\bbD_1^{-1}(z)\{\underline{m}(z)\bbT_n+{\bf I}_p\}^{-1}\bbq_i\cdot\bbq_i^{*}\bbD_1^{-1}(z)\bbq_i\\
  &=&{\bf e}_i^T\bbQ^{*}\bbD_1^{-1}(z)\{\underline{m}(z)\bbT_n+{\bf I}_p\}^{-1}\bbQ{\bf e}_i{\bf e}_i^T\bbQ^{*}\bbD_1^{-1}(z)\bbQ{\bf e}_i\\
  &=&z^{-2}{\bf e}_i^T\bbQ^{*}\{\underline{m}(z)\bbT_n+{\bf I}_p\}^{-2}\bbQ{\bf e}_i{\bf e}_i^T\bbQ^{*}\{\underline{m}(z)\bbT_n+{\bf I}_p\}^{-1}\bbQ{\bf e}_i\\
  &=&z^{-2}{\bf e}_i^T\bbQ^{*}\{\underline{m}(z)\bbT_n+{\bf I}_p\}^{-1}\bbQ\bbQ^{*}(\bbQ\bbQ^{*})^{-1}\{\underline{m}(z)\bbT_n+{\bf I}_p\}^{-1}\bbQ{\bf e}_i{\bf e}_i^T\bbQ^{*}\{\underline{m}(z)\bbT_n+{\bf I}_p\}^{-1}\bbQ{\bf e}_i\\
  &=&z^{-2}{\bf e}_i^T\bbQ^{*}\{\underline{m}(z)\bbT_n+{\bf I}_p\}^{-1}\bbQ\bbQ^{*}(\bbQ\bbQ^{*})^{-1/2}\{\underline{m}(z)\bbT_n+{\bf I}_p\}^{-1}
  (\bbQ\bbQ^{*})^{-1/2}\bbQ{\bf e}_i\\
  &&\cdot{\bf e}_i^T\bbQ^{*}\{\underline{m}(z)\bbT_n+{\bf I}_p\}^{-1}\bbQ{\bf e}_i,
\end{eqnarray*}
and
$$\bbQ^*(\bbQ\bbQ^*)^{-1/2}(\underline{m}(z)\bbQ\bbQ^{*}+\bbI_p)^{-1}(\bbQ\bbQ^*)^{-1/2}\bbQ=\bbI_k-[\bbQ^{*}\bbQ+\underline{m}^{-1}(z)\bbI_k]^{-1}\bbQ^{*}\bbQ$$
$$\bbQ^{*}(\underline{m}(z)\bbQ\bbQ^{*}+\bbI_p)^{-1}\bbQ=\underline{m}^{-1}(z)\{\bbI_k-[\bbI_k+\underline{m}(z)\bbQ^{*}\bbQ]^{-1}\}.$$
That is,
\begin{eqnarray*}
  &&\rE\umn(z)-\umn^0(z)\nonumber\\
  &=&-\frac{\umn^0(z)A_n(z)}{1-c_n\rE\umn(z)\umn^0(z)\int\frac{t^2dH_n(t)}
    {(1+t\rE\umn(z))(1+t\umn^0(z))}}\nonumber\\
  &=&
  \frac{\alpha_xy\int\frac{\underline{m}^3(z)t^2dH(t)}{(1+t\underline{m}(z))^3}}{\left(1-y\int\frac{\underline{m}^2(z)t^2dH(t)}
      {(1+t\underline{m}(z))^2}\right)\left(1-\alpha_xy\int\frac{\underline{m}^2(z)t^2dH(t)}{(1+t\underline{m}(z))^2}\right)}
  +\beta_x\frac{y\int\frac{\underline{m}^3(z)t^2}
    {(\underline{m}(z)t+1)^3}dH(t)}{1-y\int\frac{\underline{m}^2(z)t^2dH(t)}
    {(1+t\underline{m}(z))^2}}.
\end{eqnarray*}

%
\subsection{\mgai{Proof of Corollary \protect\ref{cor1}}}

Let $a(y)=(1-\sqrt{y})^2$ and $b(y)=(1+\sqrt{y})^2$, then for $\ell\in\{1,\ldots,L\}$,  we have
\begin{eqnarray*}
  EX_{f_{\ell}}&=&\frac{[a(y)]^{\ell}+[b(y)]^{\ell}}{4}-\frac{1}{2\pi}\int\limits_0^\pi(1+y-2\sqrt{y}\cos\theta)^{\ell}d\theta\\
  &&-\beta_x\frac{1}{2\pi {\bf i}}\oint z^{\ell}\frac{y\underline{m}^3(z)(1+\underline{m}(z))^{-3}}{1-y\underline{m}^2(z)(1+\underline{m}(z))^{-2}}dz\\
  &=&\frac{[a(y)]^{\ell}+[b(y)]^{\ell}}{4}-\frac{1}{2}\sum\limits_{\ell_1=0}^{\ell}
  \binom{\ell}{\ell_1}^2y^{\ell_1}+\beta_x\sum\limits_{\ell_2=2}^{\ell}\binom{\ell}{\ell_2-2}\binom{\ell}{\ell_2}y^{\ell+1-\ell_2},
\end{eqnarray*}
where
\begin{eqnarray*}
  &&\frac{1}{2\pi{\bf i}}\oint z^{\ell}\frac{y\underline{m}^3(z)(1+\underline{m}(z))^{-3}}{1-y\underline{m}^2(z)(1+\underline{m}(z))^{-2}}dz\\
  &=&\frac{y}{2\pi {\bf i}}\oint\left(-\frac{1}{\underline{m}(z)}+\frac{y}{1+\underline{m}(z)}\right)^{\ell}
  \frac{\underline{m}(z)}{(1+\underline{m}(z))^3}d\underline{m}(z)\nonumber\\
  &=&\frac{y}{2\pi{\bf i}}\oint\sum\limits_{{\ell_1}=0}^{\ell}\frac{\binom{\ell}{{\ell_1}}(-1)^jy^{\ell-{\ell_1}}\underline{m}(z)}
  {\underline{m}^{\ell_1}(z)(1+\underline{m}(z))^{\ell+3-{\ell_1}}}d\underline{m}(z)\\
  &=&\frac{y}{2\pi{\bf i}}\oint\sum\limits_{{\ell_1}=2}^{\ell}\frac{\binom{\ell}{{\ell_1}}(-1)^jy^{\ell-{\ell_1}}}
  {\underline{m}^{{\ell_1}-1}(z)(1+\underline{m}(z))^{\ell+3-{\ell_1}}}d\underline{m}(z)\\
  &=&\sum\limits_{{\ell_1}=2}^{\ell}\frac{\binom{\ell}{{\ell_1}}(-1)^jy^{\ell+1-{\ell_1}}}{2\pi {\bf i}}\oint\frac{1}{\underline{m}^{{\ell_1}-1}(z)(1+\underline{m}(z))^{\ell+3-{\ell_1}}}d\underline{m}(z)\\
  &=&-\sum\limits_{{\ell_1}=2}^{\ell}\binom{\ell}{{\ell_1}-2}\binom{\ell}{{\ell_1}}y^{\ell+1-{\ell_1}}.
\end{eqnarray*}
and
\begin{eqnarray*}
  &&\Cov(X_{f_{\ell}}, X_{f_{\ell'}})\\
  &=&2y^{\ell+\ell'}\sum\limits_{\ell_1=0}^{\ell-1}\sum\limits_{\ell_2=0}^{\ell'}
  \binom{\ell}{\ell_1}\binom{\ell'}{\ell_2}\left(\frac{1-y}{y}\right)^{\ell+\ell'}
  \sum\limits_{\ell_3=0}^{\ell-\ell_1}\ell_3\binom{2\ell-1-\ell_1-\ell_3}{\ell-1}\binom{2\ell'-1-\ell_2+\ell_3}{\ell'-1}\\
  &&+\beta_x y\sum\limits_{\ell_1=1}^{\ell}\binom{\ell}{\ell_1-1}\binom{\ell}{\ell_1}y^{\ell-\ell_1}\cdot\sum\limits_{\ell_2=1}^{\ell'}\binom{\ell'}{\ell_2-1}
  \binom{\ell'}{\ell_2}y^{\ell'-\ell_2},~~\ell, \ell'\in\{1,\ldots,L\}
\end{eqnarray*}
with
\begin{eqnarray*}
  \frac{1}{2\pi
    {\bf i}}\oint\frac{z^{\ell}}{(1+\underline{m}(z))^2}d\underline{m}(z)&=&\frac{1}{2\pi
    {\bf i}}\oint\frac{1}{(1+\underline{m}(z))^2}\left(-\frac{1}{\underline{m}(z)}+\frac{y}{1+\underline{m}(z)}\right)^{\ell}d\underline{m}(z)\\
  &=&\frac{1}{2\pi
    {\bf i}}\oint\sum\limits_{{\ell_1}=0}^{\ell}\frac{\binom{\ell}{{\ell_1}}(-1)^jy^{\ell-{\ell_1}}}
  {\underline{m}^{\ell_1}(z)(1+\underline{m}(z))^{\ell+2-{\ell_1}}}d\underline{m}(z)\\
  &=&\sum\limits_{{\ell_1}=1}^{\ell}\frac{\binom{\ell}{{\ell_1}}(-1)^jy^{\ell-{\ell_1}}}{2\pi
    i}\oint\frac{1}{\underline{m}^{\ell_1}(z)[1+\underline{m}(z)]^{\ell+2-{\ell_1}}}d\underline{m}(z)\\
  &=&\sum\limits_{{\ell_1}=1}^{\ell}\binom{\ell}{{\ell_1}-1}\binom{\ell}{{\ell_1}}y^{\ell-{\ell_1}}.
\end{eqnarray*}

\appendix
\section{Mathematical Tools}

\begin{lemma}\label{lem1} [\citet{Burkholder1973}].   Let  $\{X_k\}$ be a
  complex martingale difference sequence with respect to the increasing
  $\sigma$-field  $\{\mathscr F_k\}$.  Then for $p>1$
  $$\rE\left|\sum X_k\right|^p\leq
  K_p\rE\left(\sum|X_k|^2\right)^{p/2}.$$
\end{lemma}

\medskip

\begin{lemma}\label{lem2}
  For $\bbX=(X_1,\ldots,X_n)^T$ i.i.d. standardized (complex)
  entries, $\bbC=(c_{ij})$ $n\times n$ matrix (complex) we have
  $$\rE|\bbX^*\bbC\bbX-\rtr\bbC|^4\leq K\left(\left(
      \rtr  (\bbC\bbC^*)\right)^{2}
    +\siln \rE|X_{ii}^8||c_{ii}|^4\right). $$
\end{lemma}

The proof of the lemma can easily follow by simple calculus and thus omitted.

\begin{lemma}\label{lem3} Let $f_1,f_2,\ldots$ be analytic in $D$, a connected
  open set of ${\mathbb C}$, satisfying $|f_n(z)|\leq M$ for every $n$ and $z$ in $D$, and
  $f_n(z)$ converges, as $n\to\infty$ for each $z$ in a subset of $D$ having
  a limit point in $D$.  Then there exists a function $f$, analytic in $D$
  for which $f_n(z)\to f(z)$ and $f'_n(z)\to f'(z)$ for all $z\in D$ where $'$ denotes the derivative.
  Moreover, on any set bounded by a contour interior to $D$ the convergence is uniform
  and $\{f'_n(z)\}$ is uniformly bounded. \end{lemma}

\begin{lemma}\label{lem4} [Theorem 35.12 of \citet{Billingsley1995}]. Suppose for each
  $n$ $Y_{n1}$, $Y_{n2},\ldots,$ $Y_{nr_n}$ is a real martingale difference sequence
  with respect to the increasing $\sigma$-field  $\{\mathscr F_{nj}\}$ having
  second moments.  If as $n\to\infty$
  \be\sum_{j=1}^{r_n}\rE(Y^2_{nj}|\mathscr F_{n,j-1})\iparrow\sigma^2\label{tag i}\ee
  {\sl where $\sigma^2$ is a positive constant, and for each $\epsilon>0$}
  \be\sum_{j=1}^{r_n}\rE(Y^2_{nj}I_{(|Y_{nj}|\geq\epsilon)})\to0\label{tag ii}\ee
  {\sl then}
  $$\sum_{j=1}^{r_n}Y_{nr_n}\darrow N(0,\sigma^2).$$
\end{lemma}
\begin{lemma} \label{lem5} (Lemma 2.6 of \citet{bai1999}).  For
  $p\times n$ Hermitian $\bbA$ and $\bbB$
  $$\|F^{\bbA\bbA*}-F^{\bbB\bbB^*}\|\leq p^{-1}
  {\rm rank}(\bbA-\bbB),$$
  where $\|\cdot\|$ here denotes sup norm on functions.\medskip
\end{lemma}

\begin{lemma} \label{lem6} (Lemma 2.7 of \citet{bai1999}).  For
  $p\times n$ Hermitian $\bbA$ and $\bbB$
  $$L^4(F^{\bbA\bbA^*},F^{\bbB\bbB^*})\leq 2p^{-2}
  {\rm tr}(\bbA-\bbB)(\bbA^*-\bbB^*){\rm tr}(\bbA\bbA^*+\bbB\bbB^*),$$
  where $L(F,G)$ denotes the Levy distance between distribution functions.\medskip
\end{lemma}

\begin{lemma} \label{lem7} (Lemma 2.6 of \citet{SB95}).
  Let $z\in {\mathbb C}^+$ with $v=\Im\,z$, $\bbA$ and $\bbB$ being $n\times n$ with $\bbB$
  Hermitian, and $\bbr\in{\mathbb C} ^n$.  Then
  $$\bigl|\rtr\bigl((\bbB-z\bbI)^{-1}-(\bbB+\bbr\bbr^*-z\bbI)^{-1}\bigr)\bbA\bigr|=
  \left|\frac{\bbr^*( \bbB-z\bbI)^{-1}\bbA( \bbB-z\bbI)^{-1}\bbr}
    {1+\bbr^*(\bbB-z\bbI)^{-1}\bbr}\right|\leq\frac{\|\bbA\|}v.$$
\end{lemma}

\section{Truncation and normalization for the proof of Theorem 2.1}

\subsection{Truncation of the matrix $\bbQ$}

Because $H$ is a proper distribution function, for any given $\ep>0$, there exists a constant $\tau>0$ such that
$ 1-H(\tau)<\ep$. Without loss of generality, we may assume that $\tau$ is a continuity point of $H$. Suppose the singular value decomposition of $\bbQ$ is given by
$$
\bbQ=\bbU\bgL\bbV^*
$$
where $\bbU_{p\times p}$ and $\bbV_{k\times k}$ are two unitary matrices and $\bgL_{p\times k}=\diag[l_1,l_2,\cdots]$ is a diagonal matrix of nonnegative real singular values of $\bbQ$ in descending order. Define
$$
\widehat \bbQ=\bbU\diag[l_1\wedge \sqrt{\tau},\cdots, l_k\wedge \sqrt{\tau},\cdots]\bbV^*
$$
and $\widehat \bbB_n=\frac1n\widehat\bbQ\bbX_n\bbX_n^*\widehat\bbQ^*$.
By Lemma \ref{lem5}, we have
\bqn
\|F^{\bbB_n}-F^{\widehat\bbB_n}\|\le \frac1p\rank(\bbQ-\widehat\bbQ)\le \frac1p{\#}\{i:\, l_i^2> \tau\}
\to 1-H(\tau)<\ep.
\eqn
Therefore, we may assume that the norm of $\bbQ$ is bounded by some constant $\sqrt{\tau}$.

\subsection{Truncation}

By Assumption (a), there exists a sequence of constants $\eta_n\downarrow 0$ such that
\be
\frac1{pn\eta_n^2}\sum_{i=1}^k\sum_{j=1}^n \|\bbq_{i}\|^2\rE|x_{ij}^2|I\Big(|x_{ij}|>\eta_n \sqrt{n}/\|\bbq_i\|\Big)\to0.
\label{mmt3}
\ee
Define $\widehat x_{ij}=x_{ij}I(|x_{ij}|\le \eta_n\sqrt{n}/\|\bbq_i\|)$, $\widehat \bbX_n=(\widehat x_{ij})$ and $\widehat \bbB_n=n^{-1}\bbQ\widehat\bbX_n\widehat\bbX_n^*\bbQ^*_n$.
Applying Lemma \ref{lem5} again, we have
\begin{eqnarray*}
  \|F^{\bbB_n}-F^{\widehat \bbB_n}\|\le p^{-1}\rank(\bbX_n-\widehat\bbX_n)\le p^{-1}\sum_{i=1}^k\sum_{j=1}^n
  I(|x_{ij}|>\eta_n\sqrt{n}/\|\bbq_i\|)\to 0~~ a.s.
\end{eqnarray*}
because
\begin{eqnarray*}
  &&\rE\left(p^{-1}\sum\limits_{i=1}^k\sum\limits_{j=1}^nI(|x_{ij}|>{\eta_n\sqrt{n}}/{\|\bbq_i\|})\right)\\
  &\le& (pn\eta^2_n)^{-1}\sum\limits_{i=1}^k\sum\limits_{j=1}^n\|\bbq_i\|^2\rE|x_{ij}^2|I(|x_{ij}|>{\eta_n\sqrt{n}}/{\|\bbq_i\|})
  \to 0,
\end{eqnarray*}
\begin{eqnarray*}
  &&\rVar\left(p^{-1}\sum\limits_{i=1}^k\sum\limits_{j=1}^n
    I(|x_{ij}|>{\eta_n\sqrt{n}}/{\|\bbq_i\|})\right)\\
  &\le &p^{-2}\sum\limits_{i=1}^k\sum\limits_{j=1}^n\rP(|x_{ij}|>{\eta_n\sqrt{n}}/{\|\bbq_i\|})
  \\
  &\le&(p^2n\eta^2_n)^{-1}\sum\limits_{i=1}^k\sum\limits_{j=1}^n\|\bbq_i\|^2\rE|x_{ij}^2|I(|x_{ij}|>{\eta_n\sqrt{n}}/{\|\bbq_i\|})\\
  &=&o(p^{-1})
\end{eqnarray*}
and by Bernstein's inequality we have $P(p^{-1}\sum_{i=1}^k\sum_{j=1}^n
I(|x_{ij}|>\eta_n\sqrt{n}/\|\bbq_i\|)>\epsilon)\leq K\exp(-bp)$, for some constants $K<\infty$ and $b>0$.

\subsection{Centralization} \label{subsecA2}

Define $\widetilde \bbX_n=\widehat\bbX_n-\rE\widehat\bbX_n$ and $\widetilde\bbB_n=n^{-1}\bbQ\widetilde\bbX_n\widetilde\bbX_n^*\bbQ^*$. By Lemma \ref{lem6}, we have
\begin{eqnarray*}
  L^4(F^{\widehat\bbB_n},F^{\widetilde \bbB_n})\le 2p^{-2}n^{-2}\rtr\bbQ(\widehat \bbX_n-\widetilde\bbX_n)(\widehat \bbX_n-\widetilde\bbX_n)^*\bbQ^*
  \rtr(\bbQ\widehat\bbX_n\widehat\bbX_n^*\bbQ^*+\bbQ\widetilde\bbX_n\widetilde\bbX_n^*\bbQ^*).
\end{eqnarray*}
Notice that
\begin{eqnarray*}
  &&(pn)^{-1}\rtr(\bbQ(\widehat \bbX_n-\widetilde\bbX_n)(\widehat\bbX_n-\widetilde\bbX_n)^*\bbQ^*)\\
  &=&(pn)^{-1}\rtr\bbQ\rE\widehat\bbX_n(\rE\widehat\bbX_n)^*\bbQ^*\\
  &=&(pn)^{-1}\sum\limits_{i=1}^p\sum\limits_{j=1}^n\left|\sum_{\ell=1}^kq_{i\ell}\rE \hat x_{\ell j}\right|^2\\
  &\le&n^{-1}\sum\limits_{i=1}^p\sum\limits_{j=1}^n\sum\limits_{\ell=1}^k|q_{i\ell}|^2\sum\limits_{h=1}^k\frac{\|\bbq_h\|^2}{\eta_n^2pn}\rE|x_{hj}^2|I(|x_{hj}|>\eta_n\sqrt{n}/\|\bbq_h\|)\\
  &&(\mbox{by Cauchy-Schwarz inequality})\\
  &=&n^{-1}\sum\limits_{i=1}^p\sum\limits_{\ell=1}^k|q_{i\ell}|^2\cdot(\eta_n^2pn)^{-1}\sum\limits_{j=1}^n\sum\limits_{h=1}^k
  \|\bbq_h\|^2\rE|x_{hj}^2|I(|x_{hj}|>\eta_n\sqrt{n}/\|\bbq_h\|)\\
  &=&n^{-1}\rtr(\bbQ\bbQ^{*})(\eta_n^2pn)^{-1}\sum\limits_{j=1}^n\sum\limits_{h=1}^k\|\bbq_h\|^2\rE|x_{hj}^2|I(|x_{hj}|>\eta_n\sqrt{n}/\|\bbq_h\|)\\&=&o(1).
\end{eqnarray*}
Noticing that the above bound is non-random, to show that $L^4(F^{\widehat\bbB_n}, F^{\widetilde\bbB_n})\to 0,a.s.$, one only needs to prove that
\bqa
&&\frac1{pn}\rtr(\bbQ\widetilde\bbX_n\widetilde\bbX_n^*\bbQ^*)=\frac1{pn}\sum_{i=1}^p\sum_{j=1}^n\left|\sum_{\ell=1}^kq_{i\ell}\tilde x_{\ell j}\right|^2\non
&=&\frac1{pn}\sum_{i=1}^p\sum_{j=1}^n\sum_{\ell=1}^k|q_{i\ell}^2||\tilde x_{\ell j}^2|+\frac1{pn}\sum_{i=1}^p\sum_{j=1}^n\sum_{k_1\ne k_2}q_{ik_1}\bar q_{ik_2}\tilde x_{k_1,j}\bar{\tilde x}_{k_2,j}=O_{a.s.}(1).
\label{eq007}
\eqa
Note that
$$
\rE\left(\frac1{pn}\sum_{i=1}^p\sum_{j=1}^n\sum_{\ell=1}^k|q_{ik}^2||\tilde x_{\ell j}^2|\right)\le\frac1{pn}\sum_{i=1}^p\sum_{j=1}^n\sum_{\ell=1}^k|q_{i\ell}^2|=O(1)
$$
and
\bqa
&&\rE\left(\frac1{pn}\sum_{i=1}^p\sum_{j=1}^n\sum_{\ell=1}^k|q_{i\ell}^2|(|\tilde x_{\ell j}^2|-\rE|\tilde x_{\ell j}^2|)\right)^4\non
&\le&\frac1{p^4n^4}\sum_{j=1}^n\sum_{\ell=1}^k\left(\sum_{i=1}^p|q_{i\ell}^2|\right)^4\rE|\tilde x_{\ell j}^8|
+\frac3{p^4n^4}\left(\sum_{j=1}^n\sum_{\ell=1}^k\left(\sum_{i=1}^p|q_{i\ell}^2|\right)^2\rE|\tilde x_{\ell j}^4|\right)^2 \non
&\le&\frac{\eta_n^6}{p^4}\sum_{\ell=1}^k\|\bbq_\ell\|^2+\frac{3\eta_n^4}{p^4}\left(\sum_{\ell=1}^k\|\bbq_\ell\|^2\right)^2=o(p^{-2}).
\label{eq009}
\eqa
These inequalities simply imply $(pn)^{-1}\sum_{i=1}^p\sum_{j=1}^n\sum_{\ell=1}^k|q_{i\ell}^2||\tilde x_{\ell j}^2|=O_{a.s.}(1)$.
Furthermore,
\bqa
&&\rE\left(\frac1{pn}\sum_{i=1}^p\sum_{j=1}^n\sum_{k_1\ne k_2}q_{ik_1}\bar q_{ik_2}\tilde x_{k_1,j}\bar{\tilde x}_{k_2,j}\right)^2
\le \frac2{p^2n^2}\sum_{j=1}^n\sum_{k_1\ne k_2}\left|\sum_{i=1}^pq_{ik_1}\bar q_{ik_2}\right|^2\non
&\le&\frac2{p^2n}\rtr (\bbQ\bbQ^*)^2=O(p^{-2}),
\label{eq0010}
\eqa
which implies that $(pn)^{-1}\sum_{i=1}^p\sum_{j=1}^n\sum_{k_1\ne k_2}q_{ik_1}\bar q_{ik_2}\tilde x_{k_1,j}\bar{\tilde x}_{k_2,j}\to0,\,a.s.$
Hence the assertion (\ref{eq007}) is proved.

\subsection{Rescaling}\label{subsecA3}

Denote $\gs_{ij}^2=\rE |\tilde x_{ij}|^2$, and $y_{ij}$ are i.i.d. random variables satisfy $\rP(y_{ij}=\pm1)=\frac12$ and are independent of all $x_{ij}$'s. Define
$$
\breve x_{ij}=\begin{cases}
  \gs_{ij}^{-1}\tilde x_{ij}&\mbox{if } \gs_{ij}^2>1-\eta_n,\cr
  y_{ij}&\mbox{otherwise.}\cr
\end{cases}
$$
We further define $\breve \bbX_n=(\breve x_{ij})$ and $\breve\bbB_n=n^{-1}\bbQ\breve\bbX_n\breve\bbX_n^*\bbQ^*$.
Applying Lemma \ref{lem6} again, we have
\begin{eqnarray*}
  L^4(F^{\breve\bbB_n},F^{\widetilde\bbB_n})\le 2p^{-2}n^{-1}\rtr\bbQ(\widetilde\bbX_n-\breve\bbX_n)(\widetilde\bbX_n-\breve\bbX_n)^*\bbQ^*
  \cdot\rtr(\widetilde\bbB_n+\breve\bbB_n).
\end{eqnarray*}
We have proved in (\ref{eq007}) that $p^{-1}\rtr\widetilde\bbB_n=O_{a.s.}(1)$. Similarly, we can prove that $p^{-1}\rtr\breve\bbB_n=O_{a.s.}(1)$.
What remains is to show that
\begin{equation}
  \frac{1}{pn}\rtr\bbQ(\widetilde\bbX_n-\breve\bbX_n)(\widetilde\bbX_n-\breve\bbX_n)^*\bbQ^*
  =\frac{1}{pn}\sum_{i=1}^p\sum_{j=1}^n\left|\sum_{\ell=1}^kq_{i\ell}(\breve x_{\ell j}-\tilde x_{\ell j})\right|^2=o_{a.s.}(1).
  \label{eq008}
\end{equation}
Write $E_n=\{(i,j):\, \gs_{ij}^2<1-\eta_n\}$ and $E_{(j)}=\{i:\, (i,j)\in E_n\}$. Then
\begin{eqnarray*}
  \sum_{i=1}^p\sum_{j=1}^n\left|\sum_{\ell=1}^kq_{i\ell}(\breve x_{\ell j}-\tilde x_{\ell j})\right|^2
  &=&\sum_{i=1}^p\sum_{j=1}^n\sum_{\ell \in E_{(j)}}|q_{i\ell}^2||\breve{x}_{\ell j}-\tilde x_{\ell j}|^2\\
  &&+\sum_{i=1}^p\sum_{j=1}^n\sum_{k_1\ne k_2}q_{ik_1}\bar q_{ik_2}(\breve x_{k_1j}-\tilde x_{k_1j})
  (\bar{\breve x}_{k_2j}-\bar{\tilde x}_{k_2j}).
\end{eqnarray*}
Note that $\sum_{i=1}^p|q_{i\ell}|^2\le \|\bbq_{\ell}\|^2$ and hence
\begin{eqnarray*}
  &&\rE\left(\frac1{pn}\sum\limits_{i=1}^p\sum\limits_{j=1}^n\sum\limits_{\ell\in E_{(j)}}|q_{i\ell}^2||\breve{x}_{\ell j}-\tilde x_{\ell j}|^2\right)\\
  &\le&(pn)^{-1}\left(\sum\limits_{(i,j)\in E_n}\|\bbq_{i}\|^2\rE|y_{ij}-\tilde x_{ij}|^2+\sum_{(i,j)\not\in E_n}\|\bbq_{i}^2\|(1-\gs^{-1}_{ij})^2\right)\\
  &\le&(pn)^{-1}\left(\sum\limits_{(i,j)\in E_n}\|\bbq_{i}\|^2(1+\gs^2_{ij})+\eta_n^2n\sum_{i=1}^k \|\bbq_{i}^2\|\right)\\
  &\le&2(pn\eta^2_n)^{-1}\sum\limits_{i,j}\|\bbq_{i}\|^2(1-\gs^2_{ij})+\frac{\eta_n}{p}\rtr(\bbQ\bbQ^*) \\
  &\le&2(pn\eta^2_n)^{-2}\sum\limits_{i,j}\|\bbq_{i}\|^2\rE|x_{ij}^2|I(|x_{ij}|>\eta_n\sqrt{n}/\|\bbq_i\|)
  +\eta_n\|\bbQ\|^2\to 0.
\end{eqnarray*}
Similar to (\ref{eq009}) and (\ref{eq0010}), one can prove that
\begin{eqnarray*}
  \rE\left(\frac1{pn}\sum\limits_{i=1}^p\sum\limits_{j=1}^n\sum\limits_{\ell\in E_{(j)}}|q_{i\ell}^2|(|\breve{x}_{\ell j}-\tilde x_{\ell j}|^2-\rE|\breve{x}_{\ell j}-\tilde x_{\ell j}|^2)\right)^4
  =O(p^{-2}\eta_n^{-4})\\
  \rVar\left(\frac1{pn}\sum\limits_{i=1}^p\sum\limits_{j=1}^n\sum\limits_{k_1\ne k_2}q_{ik_1}\bar q_{ik_2}(\breve x_{k_1j}-\tilde x_{k_1j})
    (\bar{\breve x}_{k_2j}-\bar{\tilde x}_{k_2j})\right)=O(p^{-2}).
\end{eqnarray*}
From these, it is easy to show (\ref{eq008}).

\section{Truncation and normalization for the proof of Theorem 3.1}

\subsection{Truncation}

By Assumption (c), there exists a sequence of constants $\eta_n\downarrow 0$ such that
\be
\frac1{pn\eta_n^6}\sum_{i=1}^k\sum_{j=1}^n \|\bbq_{i}\|^2\rE|x_{ij}^4|I\Big(|x_{ij}|>\eta_n \sqrt{n/\|\bbq_i\|}\Big)\to0,
\label{condclt}
\ee
Define $\widehat x_{ij}=x_{ij}I(|x_{ij}|\le \eta_n\sqrt{n/\|\bbq_i\|})$, $\widehat \bbX_n=(\widehat x_{ij})$ and $\widehat \bbB_n=n^{-1}\bbQ\widehat\bbX_n\widehat\bbX_n^*\bbQ^*_n$.
Then,
\bqn
\rP(\bbB_n\ne \widehat \bbB_n)&\le&\rE\sum\limits_{i=1}^k\sum\limits_{j=1}^n
I(|x_{ij}|>\eta_n\sqrt{n/\|\bbq_i\|})\\
&\le&\eta_n^{-4}n^{-2}\sum\limits_{i=1}^k\sum\limits_{j=1}^n\|\bbq_i\|^2\rE|x_{ij}^4|I(|x_{ij}|>\eta_n\sqrt{n/\|\bbq_i\|})\to 0.
\eqn

\subsection{Centralization}

Similarly define $\widetilde \bbX_n=\widehat\bbX_n-\rE\widehat\bbX_n$ and $\widetilde\bbB_n=n^{-1}\bbQ\widetilde\bbX_n\widetilde\bbX_n^*\bbQ^*$. By
Lemma \ref{lem6}, we have
\begin{equation}
  L^4(F^{\widehat\bbB_n}, F^{\widetilde \bbB_n})\le n^{-2}p^{-2}\rtr (\bbQ(\widehat \bbX_n-\widetilde\bbX_n)(\widehat \bbX_n-\widetilde\bbX_n)^*\bbQ^*)
  \rtr(\bbQ\widehat\bbX_n\widehat\bbX_n^*\bbQ^*+\bbQ\widetilde\bbX_n\widetilde\bbX_n^*\bbQ^*).
\end{equation}
Notice that
\bqn
&&(pn)^{-1}\rtr(\bbQ(\widehat\bbX_n-\widetilde\bbX_n)(\widehat\bbX_n-\widetilde\bbX_n)^*\bbQ^*)=(pn)^{-1}\rtr\bbQ\rE\widehat\bbX_n(\rE\widehat\bbX_n)^*\bbQ^*\\
&=&(pn)^{-1}\sum\limits_{i=1}^p\sum\limits_{j=1}^n\left|\sum\limits_{\ell=1}^kq_{i\ell}\rE \hat x_{\ell j}\right|^2\quad\mbox{(by Cauchy-Schwarz inequality)}\\
&\le&(pn)^{-1}\sum\limits_{i=1}^p\sum\limits_{j=1}^n\sum\limits_{\ell=1}^k|q_{i\ell}^2|\sum\limits_{h=1}^k\frac{\|\bbq_h\|^3}{\eta_n^6n^3}
\rE|x_{\ell j}^4|I(|x_{hj}|>\eta_n\sqrt{n/\|\bbq_\ell\|})\\
&=&n^{-1}\sum\limits_{i=1}^p\sum\limits_{\ell=1}^k|q_{i\ell}^2|\cdot\sum\limits_{j=1}^n\sum\limits_{h=1}^{k}\frac{\|\bbq_h\|^3}{\eta_n^6pn^3}
\rE|x_{\ell j}^4|I(|x_{hj}|>\eta_n\sqrt{n/\|\bbq_{\ell}\|})\\
&=&o(1).
\eqn
By the similar approach given in subsection \ref{subsecA2}, one may prove that
\bqa
&&(np)^{-1}[\rtr(\bbQ\widehat\bbX_n\widehat\bbX_n^*\bbQ^*)+\rtr(\bbQ\widetilde\bbX_n\widetilde\bbX_n^*\bbQ^*)]=O_{a.s.}(1).
\label{eq008R1}
\eqa
Hence  $L^4(F^{\widehat\bbB_n}, F^{\widetilde \bbB_n})\rightarrow0$ a.s.

\subsection{Rescaling}
Denote $\gs_{ij}^2=\rE |\tilde x_{ij}|^2$, and $y_{ij}$ are i.i.d. random variables satisfy $\rP(y_{ij}=\pm1)=\frac12$. Define
$$
\breve x_{ij}=\begin{cases}
  \gs_{ij}^{-1}\tilde x_{ij}&\mbox{if } \gs_{ij}^2>1-\eta_n,\cr
  y_{ij}&\mbox{otherwise.}\cr
\end{cases}
$$
We further define $\breve \bbX_n=(\breve x_{ij})$ and $\breve\bbB_n=\frac1n\bbQ\breve\bbX_n\breve\bbX_n^*\bbQ^*$.
Similar to subsection \ref{subsecA3}, one can show that
\bqn
L^4(F^{\breve\bbB_n}, F^{\widetilde\bbB_n})\to 0, a.s.
\eqn

\bigskip\bigskip

\end{document}